\documentclass[12pt]{article}
\usepackage{amssymb,amsthm,amsmath,amscd,amsfonts,bbm,mathrsfs,fullpage,comment}
\usepackage[all]{xy}

 \def\noproof{{\unskip\nobreak\hfill\penalty50\hskip2em\hbox{}%
      \nobreak\hfill$\Box$\parfillskip=0pt%
     \finalhyphendemerits=0\par}}
\def\enddemo{\ifmmode\eqno\Box\else\noproof\vskip0.8truecm\fi}
\newtheorem{ithm}{Theorem}
\newtheorem{theo}{Theorem}[section]

\newtheorem{definition}[theo]{Definition}
\newtheorem{lemma/def}[theo]{Lemma/Definition}
\newtheorem{conjecture}[theo]{Conjecture}

\newtheorem{prop}[theo]{Proposition}
\newtheorem{corollary}[theo]{Corollary}

\newtheorem{remark}[theo]{Remark}
\newtheorem{remarks}[theo]{Remarks}

\newtheorem{lemma}[theo]{Lemma}

\newcommand{\lra}{\longrightarrow}

\DeclareMathOperator{\Eis}{Eis}

\DeclareMathOperator{\Res}{Res}
\DeclareMathOperator{\Coind}{Coind}
\DeclareMathOperator{\Ind}{Ind}

\DeclareMathOperator{\aug}{aug}

\DeclareMathOperator{\res}{res}
\DeclareMathOperator{\cor}{cor}

\DeclareMathOperator{\sgn}{\ep}

\DeclareMathOperator{\rec}{rec}

\DeclareMathOperator{\ord}{ord}

\DeclareMathOperator{\id}{id}

\DeclareMathOperator{\ev}{ev}
\DeclareMathOperator{\spez}{sp}
\DeclareMathOperator{\Hom}{Hom}

\DeclareMathOperator{\pr}{pr}

\DeclareMathOperator{\supp}{supp}

\DeclareMathOperator{\wcdot}{\, \cdot\, }

\DeclareMathOperator{\der}{d}
\DeclareMathOperator{\Ker}{Ker}

\DeclareMathOperator{\Image}{Im}

\DeclareMathOperator{\Meas}{Meas}

\DeclareMathOperator{\Real}{Re}

\DeclareMathOperator{\Sh}{Sh}

\DeclareMathOperator{\Cl}{Cl}

\DeclareMathOperator{\Gal}{Gal}
\DeclareMathOperator{\Norm}{N}

\DeclareMathOperator{\Sym}{Sym}

\DeclareMathOperator{\GL}{GL}
\DeclareMathOperator{\PGL}{PGL}

\DeclareMathOperator{\incl}{incl}
\DeclareMathOperator{\an}{an}

\DeclareMathOperator{\univ}{univ}
\DeclareMathOperator{\mult}{mult}

\DeclareMathOperator{\Pol}{Pol}
\DeclareMathOperator{\Ta}{Ta}
\DeclareMathOperator{\cycl}{cycl}

\newcommand{\Cd}{{C^{\diamond}}}
\newcommand{\barCd}{{\overline{C}^{\diamond}}}
\newcommand{\Ccd}{{C_c^{\diamond}}}
\newcommand{\cCd}{{{\mathcal C}^{\diamond}}}
\newcommand{\cCcd}{{{\mathcal C}_c^{\diamond}}}

\newcommand{\llangle}{{\langle\!\langle}}
\newcommand{\rrangle}{{\rangle\!\rangle}}

\newcommand{\fa}{{\mathfrak a}}
\newcommand{\fb}{{\mathfrak b}}
\newcommand{\fA}{{\mathfrak A}}

\newcommand{\fp}{{\mathfrak p}}
\newcommand{\fP}{{\mathfrak P}}
\newcommand{\fq}{{\mathfrak q}}

\newcommand{\bbA}{{\mathbb A}}
\newcommand{\bbI}{{\bf{\mathbb I}}}
\newcommand{\bC}{{\mathbb C}}

\newcommand{\bG}{{\mathbb G}}

\newcommand{\bN}{{\mathbb N}}

\newcommand{\bQ}{{\mathbb Q}}
\newcommand{\bR}{{\mathbb R}}

\newcommand{\bZ}{{\mathbb Z}}
\newcommand{\barQ}{{\overline{\mathbb Q}}}

\newcommand{\bA}{{\mathbb A}}

\newcommand{\bhatZ}{{\widehat{\bZ}}}

\newcommand{\barM}{{\overline{M}}}
\newcommand{\barR}{{\overline{R}}}
\newcommand{\barT}{{\overline{T}}}

\newcommand{\barX}{{\overline{X}}}

\newcommand{\barfa}{{\overline{\fa}}}

\newcommand{\barchi}{{\overline{\chi}}}
\newcommand{\barphi}{{\overline{\phi}}}
\newcommand{\barkappa}{{\overline{\kappa}}}

\newcommand{\cA}{{\mathcal A}}

\newcommand{\cC}{{\mathcal C}}

\newcommand{\cF}{{\mathcal F}}

\newcommand{\cI}{{\mathcal I}}
\newcommand{\cK}{{\mathcal K}}
\newcommand{\cL}{{\mathcal L}}

\newcommand{\cN}{{\mathcal N}}
\newcommand{\cO}{{\mathcal O}}
\newcommand{\cP}{{\mathcal P}}
\newcommand{\cQ}{{\mathcal Q}}
\newcommand{\cR}{{\mathcal R}}

\newcommand{\sE}{{\mathscr E}}
\newcommand{\sF}{{\mathscr F}}

\newcommand{\sI}{{\mathscr I}}

\newcommand{\sL}{{\mathscr L}}
\newcommand{\sM}{{\mathscr M}}

\newcommand{\sR}{{\mathscr R}}

\newcommand{\sV}{{\mathscr V}}

\newcommand{\wP}{\widetilde{P}}

\newcommand{\wchi}{\widetilde{\chi}}

\newcommand{\ep}{{\epsilon}}

\newcommand{\barpsi}{{\overline{\psi}}}

\newcommand{\bu}{{\bullet}}

\newcommand{\noi}{\noindent}

\begin{document}

\title{The Eisenstein cocycle, partial zeta values and Gross--Stark units}

\author{Samit Dasgupta and Michael Spie{\ss}}
\date{}

\maketitle

\begin{abstract}
We introduce an integral version of the Eisenstein cocycle. As applications we prove a conjecture of Gross regarding the ``order of vanishing" of Stickelberger elements relative to an abelian tower of fields and give a cohomological construction of the conjectural Gross--Stark units. 
\end{abstract}

\tableofcontents

\section{Introduction}

Let $F$ be a totally real number field of degree $n$ over $\bQ$ and let $K/F$ be a finite abelian extension with Galois group $G$. Let $S$ be the set of non places of $F$ that are ramified in $K$. For $\sigma\in G$ the partial zeta function \[ \zeta_S(\sigma,s) = \sum_{\fa} \Norm(\fa)^{-s}, \qquad s\in \bC, \Real(s)>1 \] has a meromorphic continuation to the whole complex plane. Here $\fa$ ranges through all nonzero  ideals of $\cO_F$ that are relatively prime to $S$ and are mapped to $\sigma$ under the Artin map. By a theorem of Siegel and Klingen the values $\zeta_S(\sigma,k)$ at integers $k\le 0$ are rational. Furthermore by a result of Cassou-Nogu{\`e}s and Deligne--Ribet \cite{cassou, deligneribet} a certain slight modification $\zeta_{S,T}(\sigma,s)$ of $\zeta_S(\sigma,s)$ has integral values at non-positive integers $s=k$ (here $T$ is a finite set of prime numbers satisfying certain properties; see section \ref{subsection:totalreal} for the definition of $\zeta_{S,T}(\sigma,s)$ and further details). 

Consider the Stickelberger element \[ \zeta_{S, T}(K/F, k) = \sum_{\sigma\in G}\zeta_{S,T}(\sigma, k)[\sigma^{-1}]\in \bZ[G]. \]
Our main result concerns the ``order of vanishing" of  
$\zeta_{S, T}(K/F, 0)$ with respect to an intermediate field $E$ of $K/F$. Let $I$ denote the kernel of the projection $\bZ[G]\to \bZ[\Gal(E/F)]$. Let $S_{\infty}$ be the set of archimedean places of $F$ and let $r$ denote the number of places in $S\cup S_{\infty}$ that split completely in $E$. Our first main result (see Cor.\ \ref{corollary:grosstower}) is:
\begin{ithm} We have \label{mainresult}
\[ \zeta_{S, T}(K/F, 0)\, \in\, I^r \]
unless $E$ is totally real, in which case we  have:  
\[ \zeta_{S, T}(K/F, 0)\, \in\, I^{r-1} \qquad \text{ and } \qquad 2\zeta_{S, T}(K/F, 0)\, \in\, I^{r}. \]
\end{ithm}
This settles a conjecture of Gross except in the case that $E$ is totally real, where
Gross predicts the slightly stronger statement $\zeta_{S, T}(K/F, 0)\, \in\, I^{\min\{r,\sharp S \cup S_\infty-1\}}$
(see the remarks following Cor.\ \ref{corollary:grosstower}).
 In fact, in Theorem \ref{theo:highervanishing}(a) below we obtain a strengthening of Theorem~\ref{mainresult} and also obtain some related statements for $\zeta_{S, T}(K/F, k)$ when $k$ is negative.
 
Theorem~\ref{mainresult} is proven through a study of an integral version of the Eisenstein cocycle. The Eisenstein cocycle is a certain element of the $(n-1)$-st group cohomology of $\GL_n(\bQ)$. Variants of it as well as its relation to partial zeta-values have been studied by several authors \cite{stevens, sczech, nori, hill, spiess2, cdg}. Our framework differs slightly from e.g.\ \cite{cdg}, \cite{spiess2}. Let $n\ge 2$ be an integer and let $T$ be a finite set of prime numbers containing at least two elements. Denote by $\Gamma^T$ the subgroup of $\GL_n(\bQ)$ that 
stabilizes the lattices $\bZ_p^n$ and $p\bZ_p\oplus \bZ_p^{n-1}$ for all $p\in T$. In section \ref{section:eisenstein} we define the {\it integral Eisenstein cocycle}
\[ 
\Eis\,\in H^{n-1}(\Gamma^T,\sM(\cP_c(\bQ^n, \bZ)^{\sL}, \bZ)(\sgn)).
\]
The coefficients consists roughly of the dual of $\bZ$-valued locally polynomial functions on $\bA_{\bQ}^n$ with compact support, where $\bA_{\bQ}$ is the adele ring of $\bQ$ (see section \ref{subsection:cassou} for the precise definition of $\sM(\cP_c(\bQ^n, \bZ)^{\sL}, \bZ)(\sgn)$). 

By taking the cap-product of the Eisenstein cocycle with certain homology classes in degree $n-1$ one recovers the values of partial zeta functions of totally real fields at non-positive integers. To state this precisely, let $\bA=\bA_F$ denote the ring of adeles of $F$ and let $U^S$ denote the subgroup of $\bA^*$ of ideles $(x_v)_v$ with local components $x_v=1$ if $v\in S$, $x_v>0$ if $v\mid\infty$ and $x_v$  a local unit if $v\not\in S\cup S_{\infty}$. We construct natural pairings 
\begin{eqnarray}
\label{capdelta}
&& H^{n-1}(\Gamma^T,\sM(\cP_c(\bQ^n, \bZ)^{\sL}, \bZ)(\sgn)) \times H_{n-1}(F^*, C_c^0(\bA^*/U^S, \bZ[G])) \stackrel{\cap_k}{\lra} \,\,  \bZ[G] 
\end{eqnarray}
for all $k\in \bZ_{\le 0}$, where $C_c^0$ denotes the space of locally constant functions with compact support. Moreover there exists a canonical element $\rho_{K/F}\in H_{n-1}(F^*, C_c^0(\bA^*/U^S, \bZ[G]))$ defined in terms of the reciprocity map of class field theory for the extension $K/F$. We then show (see Prop.\ \ref{prop:zetaelement}) that $\zeta_{S, T}(K/F, k)= \Eis\cap_k \rho_{K/F}$ for all $k\in \bZ_{\le 0}$.

The proof of Theorem \ref{mainresult} is based on the key vanishing results Prop.\ \ref{prop:higherordervanishing} and Prop.\ \ref{prop:higherordervanishing2} in section \ref{section:chap3}. To explain our method in more detail assume for simplicity that $E$ is totally imaginary and that $S_1\subseteq S$ is the set of places in $S$ that split completely in $E$. Since $\bA^*\subseteq \prod_{v\in S_1} F_v \times (\bA^{S_1})^*$ we obtain a map 
\[
\iota: C_c^0(\bA^*, \bZ[G]) \, \to\, C_c^0(\prod_{v\in S_1} F_v\times (\bA^{S_1})^*, \bZ[G]/I^r)
\]
given by first extending a locally constant function with compact support $f\colon \bA^*\to \bZ[G]$ by zero $f_!:\prod_{v\in S_1} F_v\times  (\bA^{S_1})^*\to \bZ[G]$ and then reducing its values modulo $I^r$. Prop.\ \ref{prop:higherordervanishing} implies that $\rho_{K/F}$ lies in the kernel of the induced map 
\[
\iota_*:H_{n-1}(F^*, C_c^0(\bA^*, \bZ[G])) \to H_{n-1}(F^*, C_c^0(\prod_{v\in S_1} F_v\times (\bA^{S_1})^*, \bZ[G]/I^r)).
\]
Moreover it is not hard to see that the map 
\[
 H_{n-1}(F^*, C_c^0(\bA^*/U^S, \bZ[G])) \to \bZ[G]/I^r, \qquad x\mapsto \Eis\cap_0 x \mod I^r
 \]
factors through $\iota_*$.   Combining these facts gives the desired result $\zeta_{S, T}(K/F, 0) \in I^r$.

Prop.\ 3.6 is a refinement of Prop.\ 4.6 in \cite{spiess1}. Our proof of Prop.\ 3.6 is simpler and conceptually more satisfactory than the proof in loc.\ cit.\ as we are able to avoid the combinatorial arguments given there altogether.

Our next main application of the Eisenstein cocycle is a 
conjectural construction of  Gross--Stark units related to the ``leading term" of the 
Stickelberger element $\zeta_{S, T}(K/F, 0)$.
To explain this construction let $K/F$ be as above and assume that $\fp$ is a nonarchimedean place of $F$ that splits completely in $K$. 
Let $F_{\fp}$ denote the completion of $F$ at $\fp$ and let $\cO_{\fp}$ be its valuation ring. We define $u= u_{K/F,\fp}\in K_{\fp} = K\otimes_F F_{\fp}$ by
\[
u= \{u_\sigma\}_{\sigma\in G} \, =\, \Eis \cap (c\cap \rho^{\fp}_{K/F})\in \prod_{\sigma\in G} F_{\fp} \cong \prod_{\fP\mid\fp} K_{\fP}\cong K_{\fp}.
\]
Here the class $\rho^{\fp}_{K/F}\in H_n(F^*, C_c^0((\bA^{\{\fp\}})^*/U^{\{\fp\}\cup S}, \bZ[G]))$ is defined similarly to $\rho_{K/F}$ in terms of the reciprocity map for the extension $K/F$. The class $c=c_{\id}\in H^1(F^*, C_c(F_{\fp}, F_{\fp}^*))$ is the class of the $1$-cocycle $x\mapsto (1-x)(1_{\cO_{\fp}}\cdot \id_{F_{\fp}^*})$ (see \eqref{cycle} for details). The first cap-product is a variant of the pairing \eqref{capdelta} above for $k=0$ (however it should be mentioned that contrary to \eqref{capdelta} it involves a form of $p$-adic integration). We conjecture (see Conj.\ \ref{conjecture:grossstark}) that $u$ is a global $\fp$-unit of $K$, i.e.\ that it lies in $K\subseteq K_{\fp}$ and is a unit at all finite places $v\nmid \fp$. Proposition \ref{prop:brumerstark}(a) and (e) imply that this conjecture is a further refinement of Gross' strengthening of the Brumer--Stark conjecture (\cite{gross2}, conj.\ 7.6). The element $u$ has been previously constructed in \cite{dasgupta1} in a complicated but rather explicit form in terms of values of Shintani zeta functions at $s=0$. The construction given here is simpler, conceptually more satisfactory and the proofs of the main properties Prop.\ \ref{prop:brumerstark}(a), (b), (e) of $u$ are easier than those of the corresponding statements in loc.\ cit.\ Furthermore we show that  $u$ is trivial if $K$ has a real archimedean place (see Prop.\ \ref{prop:brumerstark}(d)), a fact which was not known before and which provides further evidence for the validity of Conj.\ \ref{conjecture:grossstark}.  Note
that Prop.\ \ref{prop:brumerstark}(e) again rests on Prop.\ 3.6.

\begin{remark} \rm (a)  Burns has proven in \cite[Corollary 4.1]{burns} that Gross's conjecture follows from the Equivariant Tamagawa Number Conjecture (ETNC).  As a result, for $F=\bQ$ one deduces the full conjecture including results relating the image of $ \zeta_{S, T}(K/\bQ, 0)$ in $I^r/I^{r+1}$ to circular units, since ETNC for Dirichlet motives is known over $\bQ$.

\noi (b) Assume that $K$ is totally complex (i.e.\ $s=0$) and that $S\supseteq S_p$ where $p$ is an odd prime. Greither and Popescu \cite{popescu} have previously shown that the $p$-part of Gross' conjecture holds,
i.e.\  $\Theta_{S,T}(K/F, 0)\in I^r\otimes \bZ_p$. 

\noi (c) An alternate, more geometric construction of the Eisenstein cocycle has been recently given by Beilinson, Kings, and Levin \cite{bkl}.  
It is interesting to consider whether the methods of this paper can be combined with their construction to yield further results.
\end{remark}

{\it Acknowledgements.} We thank Cristian Popescu for helpful communications and in particular for his suggestion to strengthen our main result in the form Theorem \ref{theo:highervanishing}(a). M.S. would like to thank  the Institute for Advanced Study for hospitality.

\paragraph{\bf Notation.}

Throughout this paper all rings are commutative and have an identity element unless specified otherwise.
Given a ring $R$ and an abelian group $\Lambda$ we denote by $R[\Lambda]$ the group ring of $\Lambda$ with coefficients in $R$. For $\lambda\in \Lambda$ we denote by $[\lambda]$ the corresponding element in $R[\Lambda]$. We denote by $\ep: R[\Lambda]\to R, [\lambda]\mapsto 1$ the augmentation map and by $R[\![\Lambda]\!]$ the completion of $R[\Lambda]$ with respect to the augmentation ideal \[ I=I(\Lambda) = I_R(\Lambda)= \Ker(\aug: R[\Lambda]\to R),\] i.e.\ 
\[
R[\![\Lambda]\!] = \underset{\underset{n}{\longleftarrow}}{\lim} \, R[\Lambda]/I^n.
\]
By abuse of notation the augmentation ideal of $R[\![\Lambda]\!]$ i.e.\ the completion $\underset{\underset{n}{\longleftarrow}}{\lim} \, I/I^n$ will be denoted by $I(\Lambda)$ as well.

For a group $G$ a subgroup $H$ there exists morphisms of $\delta$-functors
\begin{equation*}
\res^G_H: H^{\bu}(G, \wcdot) \lra H^{\bu}(H, \wcdot), \quad \cor^G_H: H_{\bu}(H, \wcdot) \lra H_{\bu}(G, \wcdot)
\end{equation*}
which in degree $0$  for a $G$-module $M$ are the canonical inclusion $M^G\hookrightarrow M^H$ and projection $M_H \to M_G$ respectively.

Throughout the paper $F$ denotes a totally real number field of degree $n$ over $\bQ$ with ring of integers $\cO_F$. Let $E_F = \cO_F^*$ denote the group of global units. More generally for a finite set $S$ of nonarchimedean places of $F$ we denote by $E_S = E_{F,S}$ the group of $S$-units of $F$. For a non-zero ideal $\fa\subseteq \cO_F$ we set $N(\fa) = \sharp(\cO_F/\fa)$. For a prime number $q$ we shall write $S_q = S_q(F)$ for the set of places above $q$. By $S_{\infty}$ we denote the set of archimedean places of $F$. We denote by $\sigma_1, \ldots, \sigma_d$ the different embeddings of $F$ into $\bR$.  Places of $F$ will be denoted by $v,w$ or also by $\fp, \fq$ if they are finite. If $\fp$ is a finite place of $F$ we denote the corresponding prime ideal of $\cO_F$ also by $\fp$. For a place $v$ of $F$ we denote by $F_v$ the completion of $F$ at $v$. If $v$ is finite then $\cO_v$ denotes the valuation ring of $F_v$ and $\ord_v$ the corresponding the normalized (additive) valuation on $F_v$ (so $\ord_v(\varpi) =1$ if $\varpi\in \cO_v$ is a local uniformizer at $v$). Also we let $|\wcdot|_v$ be the associated normalized multiplicative valuation on $F_v$. Thus if $v\in S_{\infty}$ corresponds to the embedding $\sigma\colon F \to \bR$ then $|x|_v = |\sigma(x)|$ and if $v=\fq$ is finite then $|x|_{\fq} = N(\fq)^{-\ord_{\fq}(x)}$. For a place $v$ of $F$ we put $U_v = \bR_+ = \{x\in \bR\mid x>0\}$ if $v\mid \infty$ and $U_v = \cO_v^*$ if $v$ is finite. 

Let $\bbA= \bbA_F$ be the adele ring of $F$. For a set $S$ of places of $F$ we let $\bbA\!^S$ denote the $S$-adeles. We also define $U^S =\prod_{v\not\in S} U_v$,  and $U_S = \prod_{v\in S} U_v$. If $S$ contains all archimedean places then the factor group $(\bbA^S)^*/U^S$ is canonically isomorphic to the group $\cI^S$ of fractional $\cO_F$-ideals that are relatively prime to all places in $S$. We sometimes view $F$ as a subring of $\bbA\!^S$ via the diagonal embedding. 

If $T$ is a subset of $\{ 2,3,5, \ldots, \infty\}$ and $S$ is the set of places of $F$ which lie above an element of $T$ then we often write $\bbA_T$, $\bbA\!^T$ etc.\ for $\bbA_S$, $\bbA\!^S$ etc. We also write $U^q$, $U_q$, $U^{q, S}$, $U^{q,\infty}$ etc.\ for $U^{\{q\}}$, $U_{\{q\}}$, $U^{S_q\cup S}$, $U^{S_q\cup S_{\infty}}$ etc. and use a similar notation for adeles. Thus for example $\bbA^{S,\infty}$ denotes the ring $S\cup S_{\infty}$-adeles of $F$ where $S$ is a set of nonarchimedean places of $F$. For $\ell\in\{ 2,3,5, \ldots, \infty\}$ we also put $F_{\ell}= F\otimes \bQ_{\ell} = \prod_{v\in S_{\ell}} F_v$. More generally for an abelian group $A$ we put $A_{\ell}=A\otimes\bQ_{\ell}$. We shall denote by $F^*_+$, $E_+= E_{F,+}$, $E_{S,+}$ etc.\ the totally positive elements in $F$, $E_F$, $E_{S}$ etc.\

\section{Functions and measures on totally disconnected locally compact spaces}
\label{section:meas}

If $X$ and $Y$ are topological spaces then $C(X,Y)$ denotes the set of continuous maps $X\to Y$. If $R$ is a topological ring we let 
$C_c(X,R)$ denote the subset of  $C(X,R)$ of continuous maps with compact support. If we consider $Y$ (resp.\ $R$) with the discrete topology then we shall also write $C^0(X,Y)$ (resp.\ $C^0_c(X,R)$) instead of $C(X,Y)$ (resp.\ $C_c(X,R)$).

Assume now that $X$ is a totally disconnected topological Hausdorff space and $A$ a locally profinite group (the examples we have in mind for $X$ are the additive or multiplicative group of a nonarchimedean local field, or the group of ($S$-)ideles or adeles or spaces of the form $\prod_{v\in S} F_v\times \prod_{v\not\in T}' F_v^*$ where $S\subset T$ are finite sets of places of a number field $F$; the group $A$ is typically either discrete or the multiplicative or additive group of a finite $\bQ_p$-algebra). We define subgroups $\Cd(X, A) \subseteq C(X, A)$ and $\Ccd(X, A)\subseteq  C_c(X, A)$ by
\begin{eqnarray*}
&& \Cd(X, A)=C^0(X, A)+ \sum_K C(X, K), \\
&& \Ccd(X, A)=C_c^0(X, A)+ \sum_K C_c(X, K)
\end{eqnarray*}
where the sums are taken  over all compact open subgroups $K$ of $A$. So $\Ccd(X, A)$  is the subgroup of $C_c(X, A)$ generated by locally constant maps with compact support $X\to A$ and by continuous maps with compact support $X\to K\subseteq A$ for some compact open subgroup $K\subseteq A$. Similarly, $\Cd(X, A)$ is the subgroup of $C(X, A)$ generated by locally constant maps $X\to A$ and by continuous maps $X\to K\subseteq A$ for some $K$. 

If $A$ discrete then obviously $\Cd(X, A)= C(X,A)$ and $\Ccd(X, A)= C_c(X,A)$. Also if $A$ is the additive group of a finite $\bQ_p$-algebra then $A$ can be written as a union of compact open subgroups hence we have $\Ccd(X, A) = C_c(X, A)$. On the other hand if $A=\bQ_p^*$ we have  $\Cd(X, A) \ne C_c(X, A)$ in general. 

\begin{remark}
\label{remark:charcd}
\rm If $\fA$ is either the multiplicative group of a nonarchimedean local field or the group of $S$-ideles of a number field $F$ (where $S$ is a finite set of places) and if $\chi\colon \fA\to A$ is continuous homomorphism then we have $\chi\in \Cd(\fA, A)$.
\end{remark}

If $f:X\to Y$ is a continuous map between totally disconnected topological Hausdorff space are then $C(Y, A)\to C(X, A), g\mapsto g\circ f$ maps $\Cd(Y, A)$ into $\Cd(X, A)$ i.e.\ $f$ induces a homomorphism 
\begin{equation}
\label{pullback}
f^*:\Cd(Y, A)\, \lra \, \Cd(X, A), \,\, f\mapsto f^*(g) = g\circ f.
\end{equation}
Also if $U\subseteq X$ is an open subset then $C_c(U, A)\to C_c(X, A), f\mapsto f_!$ (where $f_!:X\to A$ is the extension by zero of $f\colon U\to A$) maps $\Ccd(U, A)$ into $\Ccd(X, A)$ so we get a monomorphism
\begin{equation}
\label{extensionsbyzero}
\Ccd(U, A)\, \lra \, \Ccd(X, A), \,\, f\mapsto f_!.
\end{equation}
If $A_1, A_2, A$ are locally profinite groups and $\beta\colon A_1\times A_2\to A,(a_1,a_2) \mapsto \beta(a_1, a_2)$ is a continuous bilinear map then it induces a pairing 
\begin{equation}
\label{multfunctions}
C(X, A_1)\times C(X, A_2) \, \lra \, C(X, A), \,\, (f,g) \mapsto \beta(f,g)
\end{equation}
given by $\beta(f,g)(x) = \beta(f(x), g(x))$ for all $x\in X$. It is easy to see that for $f\in \Cd(X, A_1)$ (resp.\ $f\in \Ccd(X, A_1)$) and $g\in \Cd(X, A_2)$ 
we have $\beta(f,g)\in \Cd(X, A)$ (resp.\ $\beta(f,g)\in \Ccd(X, A)$), i.e.\ \eqref{multfunctions} induces a pairings 
\begin{eqnarray}
\label{multfunctions2}
&& \Ccd(X, A_1)\times \Cd(X, A_2) \, \lra \, \Ccd(X, A),\\
\label{multfunctions2a}
&& \Cd(X, A_1)\times \Cd(X, A_2) \, \lra \, \Cd(X, A).
\end{eqnarray}
In particular if $X = X_1\times X_2$ where $X_1$, $X_2$ totally disconnected topological Hausdorff spaces then $\beta$ induces a map
\begin{equation}
\label{multfunctions3}
\Cd(X_1, A_1)\otimes \Cd(X_2, A_2) \, \lra \, \Cd(X_1\times X_2, A), \,\, f\otimes g \mapsto f\odot_{\beta} g
\end{equation}
where $(f\odot_{\beta} g)(x_1,x_2) = \beta(f(x_1), g(x_2))$ for all $(x_1,x_2)\in X_1\times X_2$. It is easy to see that \eqref{multfunctions3} maps $\Ccd(X_1, A_1)\otimes \Ccd(X_2, A_2)$ into $\Ccd(X_1\times X_2, A)$ so we have a pairing 
\begin{equation}
\label{multfunctions4}
\Ccd(X_1, A_1)\otimes \Ccd(X_2, A_2) \, \lra \, \Ccd(X_1\times X_2, A).
\end{equation}
We will consider particularly the case where $X_1$ is discrete $A_1 = \bZ$, $A_2 = A$ and $\beta: \bZ\times A\to A$ is the map $\beta(n,a) = n\cdot a$. In this case it is easy to see that \eqref{multfunctions4} is an isomorphism 
\begin{equation}
\label{multfunctions5}
C_c(X_1, \bZ)\otimes \Ccd(X_2, A) \, \lra \, \Ccd(X_1\times X_2, A).
\end{equation}

Next attach to a homomorphism $\mu\colon C_c(X, \bZ)\to \bZ$ an $A$-valued measure on $X$ for any abelian group $A$. Firstly, by tensoring $\mu$ with the identity we obtain a homomorphism $\mu_G: C_c(X, \bZ)\otimes G = C_c^0(X, G)\to G$ for any abelian group $G$. Thus if $A$ is profinite we can consider the homomorphism 
\[
\projlim_{K} \mu_{A/K} :  \projlim_{K} C_c(X, A/K)\to \projlim_{K} A/K\, =\, A
\]
where $K$ ranges over the open subgroups of $A$.  Since $C_c(X, A) \subseteq  \projlim_{K} C_c(X, A/K)$ we see that $\mu_A$ extends canonically to a homomorphism $C_c(X, A)\to A$ (which we denote by $\mu_A$ as well). For a general $A$ (not necessarily profinite) we have seen that $\mu$ induces a homomorphism $C_c(X, K)\to K$  for every compact open subgroup $K \subset A$ and we still have a homomorphism $C_c^0(X, A)\to A$. Combining these maps we see that $\mu$ induces a canonical homomorphism $\mu_A: \Cd(X, A)\, \lra\, A$. We put $\Meas(X, A) = \Hom(\Ccd(X, A),A)$. An element of $\Meas(X, A)$ is called an {\it $A$-valued measure} on $X$. Thus we obtain a homomorphism
\begin{equation}
\label{ZtoAvaluedmeasure}
\Hom(C_c(X, \bZ), \bZ)\, \lra \, \Meas(X, A),\,\, \mu \mapsto \mu_A.
\end{equation}

\section{Homology and cohomology classes associated with Hecke characters}\label{section:chap3}

Let $F$ be a totally real number field of degree $n$ over $\bQ$. We chose an ordering of the characters $\sigma_1, \ldots, \sigma_n$ of the embeddings $F\hookrightarrow \bR$. Given two arbitrary finite, disjoint sets $\Sigma_1, \Sigma_2$ of places of $F$ and a locally profinite group $A$ we put 
\begin{equation}
\label{adelicfunctions}
\cC_{?}(\Sigma_1, A)^{\Sigma_2} \, = \, C_{?}( (\bA_F^{\Sigma_2})^*/ U^{\Sigma_1\cup \Sigma_2}, A)
\end{equation}
for $?\in \{\emptyset, \diamond, c\}$. Also we sometimes write $\cC_{?}^0 (\Sigma_1, A)^{\Sigma_2}$ if $A$ is equipped with the discrete topology.
These groups carry an $(\bA_F^{\Sigma_2})^*$-action given by $(\alpha \phi)(x) = \phi(\alpha^{-1}x)$. As a special case of the map \eqref{multfunctions2} of the previous section we have a $(\bA_F^{\Sigma_2})^*$-equivariant pairing 
\begin{equation}
\label{delta1}
\cC_c(\emptyset, \bZ)^{\Sigma_2}\times \cC(\Sigma_1, A)^{\Sigma_2}\, \lra \, \cC(\Sigma_1, A)^{\Sigma_2}, \,\, (\phi, \psi) \mapsto \phi\cdot \psi
\end{equation}
where $\phi\cdot \psi$ denotes the function $xU^{\Sigma_1\cup \Sigma_2}\mapsto \phi(xU^{\Sigma_1\cup \Sigma_2}) \psi(xU^{\Sigma_2})$. 

\subsection{The  classes $\vartheta$ and $\vartheta^S$}  Next we construct a canonical class $\vartheta\in H_{n-1}(F^*, \cC_c(\emptyset, \bZ))$. By Dirichlet's unit theorem the homology group $H_{n-1}(E_+, \bZ)$ is free of rank one. Due to the chosen ordering of the embeddings $F\hookrightarrow \bR$ there is a canonical choice of  generator $\eta$ of $H_{n-1}(E_+, \bZ)$ (see e.g.\ \cite{spiess2}). Let $\cF\subseteq \bA_F^*/U$ be a fundamental domain for the action of $F^*/E_+$ on $\bA_F^*/U$. We have 
\[ \cC_c(\emptyset, \bZ) = C_c(\bA_F^*/U, \bZ)= \Ind^{\bA^*}_U \bZ \cong \Ind^{F^*}_{E_+} C(\cF,\bZ) \] as $F^*$-modules,  hence 
\begin{equation*}
H_{n-1}(E_+, C(\cF,\bZ)) \, \cong \, H_{n-1}(F^*, \cC_c(\emptyset, \bZ))
\end{equation*}
by Shapiro's Lemma.  We define $\vartheta$ corresponding to $1_{\cF}\cap\eta\in H_{n-1}(E_+, C(\cF,\bZ))$ under this isomorphism. It is easy to see that we have $x\vartheta= \vartheta$ for all $x\in \bA_F^*$.

More generally let $S$ be a finite set of nonarchimedean places of $F$ of cardinality $r$. The group $E_{S,+}$ is free-abelian of rank $n+r-1$, hence $H_{n+r-1}(E_{S,+}, \bZ)\cong \bZ$. We choose an ordering $\fp_1, \ldots, \fp_r$ of the primes in $S$ which provides us with a distinguished generator $\eta_{S}$ of $H_{n+r-1}(E_{S,+}, \bZ)\cong \bZ$. Let $\cF$ be a fundamental domain for the action of $F^*/E_{S,+}$ on $(\bA^{S}_F)^*/U^{S}$.
We define $\vartheta^{S}\in H_{n+r-1}(F^*, \cC_c(\emptyset, \bZ)^{S})$ as the homology class corresponding to $1_{\cF}\cap\eta_{S}$ under the isomorphism
\begin{equation*}
H_{n +r-1}(E_{S,+}, C(\cF,\bZ)) \, \cong \, H_{n + r-1}(F^*, C_c((\bA^{S}_F)^*/U^{S}, \bZ))
\end{equation*}
that is induced by $C_c((\bA^{S}_F)^*/U^{S}, \bZ) \cong \Ind^{F^*}_{E_{S,+}} C(\cF,\bZ)$. 
 
We want to establish the relation between $\vartheta$ and $\vartheta^S$. For that we construct a certain cohomology class $c_{\fp}\in H^1(F^*, C_c(F_{\fp}^*/U_{\fp}, \bZ))$ associated to a given nonarchimedean place $\fp$ of $F$. Consider the short exact sequence of $F_{\fp}^*$-modules 
\begin{equation}
\label{cycle0}
\begin{CD}
0@>>> C_c(F_{\fp}^*,\bZ) @> \eqref{extensionsbyzero} >> C_c(F_{\fp},\bZ) @> g\mapsto g(0) >> \bZ @>>>0.
\end{CD}
\end{equation}  
Upon taking $U_\fp$-invariants, this sequence remains exact.  Using the associated long exact sequence in $F^*$-cohomology,
we define $c_{\fp}\in H^1(F^*, C_c(F_{\fp}^*/U_{\fp}, \bZ))$ to be the image of $1\in \bZ$ under the connecting homomorphism $\bZ\to H^1(F^*, C_c(F_{\fp}^*/U_{\fp},\bZ))$. Thus $c_{\fp}$
is the cohomology class of the $1$-cocycle $z_{\fp}\colon F^*\to C_c(F_{\fp}^*/U_{\fp}, \bZ)$ given by
\begin{equation*}
z_{\fp}(x)(y) \,=\, ((1-x)1_{\cO_{\fp}})(y) \, =\, 1_{\cO_{\fp}}(y)-1_{x\cO_{\fp}}(y)
\end{equation*}
for $x\in F^*$ and $y\in F_{\fp}$. The proof of the following lemma is straightforward and will be left to the reader.

\begin{lemma} 
\label{lemma:changeS}
For $S'=\{\fp_2, \ldots, \fp_r\}$ we have $\vartheta^{S'}\,\,=\, \, c_{\fp_1}\cap \vartheta^{S}$. In particular we have \[ \vartheta\,\,=\, \, (c_{\fp_1}\cup\ldots \cup c_{\fp_1})\cap \vartheta^{S}. \]
\end{lemma}

Here the cap-product is induced by the canonical map \[ C_c(F_{\fp_1}^*/U_{\fp_1}, \bZ) \times C_c((\bA^{S}_F)^*/U^{S}, \bZ) \to C_c((\bA^{S'}_F)^*/U^{S'}, \bZ) \] (compare \eqref{multfunctions5}).

\subsection{1-cocycles attached to homomorphisms}  
Let $\psi\colon F_{\fp}^* \to A$ be continuous homomorphism where $\fp$ is a nonarchimedean place of $F$ and $A$ is a locally profinite group (the group law of $A$ will be written additively).
Note that $\psi\in \Cd(F_{\fp}^*,A)\subseteq C(F_{\fp}, A)$.  We generalize the construction of $c_\fp$ above by attaching a cohomology class
$c_\psi \in H^1(F_\fp^*, \Cd(F_\fp, A))$ to the homomorphism $\psi$.
We let $c_{\psi}$ be the class of the cocycle $z_{\psi}: F_{\fp}^*\to \Cd(F_{\fp}, A)$ defined in \cite{spiess1}, i.e.\ $z_{\psi}(x) = (1-x) (1_{\cO_{\fp}}\cdot \psi)$, or more precisely
\begin{equation}
\label{cycle}
z_{\psi}(x)(y) \,= \, 1_{x\cO_{\fp}}(y)\cdot \psi(x) + ((1_{\cO_{\fp}}- 1_{x\cO_{\fp}})\cdot \psi)_!(y) 
\end{equation}
for $x\in F_{\fp}^*$ and $y\in F_{\fp}$. By abuse of notation for a subgroup $H$ of $F_{\fp}^*$ we shall write $c_{\psi}$ also for $\res^{F_{\fp}^*}_H(c_{\psi})\in H^1(H, \Cd(F_{\fp}, A))$.

\begin{remark} 
\label{remark:classord}
\rm If $\psi=\ord_{\fp}\colon F_{\fp}^* \to \bZ$ is the normalized valuation on $F_{\fp}$, the class $c_{\ord_{\fp}}\in H^1(F^*, C_c(F_{\fp}, \bZ))$ 
is the image of $c_{\fp}$ under the homomorphism \[ H^1(F_{\fp}^*, C_c(F_{\fp}^*/U_{\fp}, \bZ)) \to H^1(F_{\fp}^*, C_c(F_{\fp}, \bZ))\] induced by the $F_{\fp}^*$-equivariant map
\begin{equation*}
C_c(F_{\fp}^*/U_{\fp}, \bZ) \lra C_c(F_{\fp}, \bZ), \qquad 1_{xU_{\fp}}\mapsto 1_{xO_{\fp}}.
\end{equation*}
\end{remark}

Let us now give a more conceptual definition of the cohomology class $c_\psi$ attached to a homomorphism $\psi\colon F_\fp^* \to A$.
Consider the short exact sequence of $F_{\fp}^*$-modules 
\begin{equation}
\label{cycle2}
\begin{CD}
0@>>> \Ccd(F_{\fp},A)@>\beta >> \widetilde{E} @>\gamma >> \Cd(F_{\fp}^*, A) @>>>0.
\end{CD}
\end{equation}
given by the push-out of the sequence \eqref{cycle0}---tensored by $\Cd(F_{\fp}^*, A)$---with respect to the homomorphism
\begin{equation}
\label{cycle3}
\begin{CD}
C_c(F_{\fp}^*, \bZ)\otimes \Cd(F_{\fp}^*, A) @> \eqref{multfunctions2} >> \Ccd(F_{\fp}^*, A) @> \eqref{extensionsbyzero} >> \Ccd(F_{\fp}, A) 
\end{CD}
\end{equation}
(as usual we define a $F_{\fp}^*$-action on $C(F_{\fp}, A)$ by $(xf)(y) = f(x^{-1}y)$).  We have a commutative diagram 
\begin{equation*}
\begin{CD}
C_c(F_{\fp}^*,\bZ)\otimes \Cd(F_{\fp}^*, A)  @>\eqref{extensionsbyzero}\otimes \id >> C_c(F_{\fp},\bZ) \otimes \Cd(F_{\fp}^*, A)\\
@VV \eqref{cycle3} V @VV \alpha V\\
\Ccd(F_{\fp},A)@> \beta >> \widetilde{E} 
\end{CD}
\end{equation*}
where the horizontal maps are injective and their cokernels are isomorphic to $\Cd(F_{\fp}^*, A)$. We identify $A$ with the $F_{\fp}^*$-submodule of $\Cd(F_{\fp}^*, A)$ of constant maps $F_{\fp}^* \to A$ and put $\barCd(F_{\fp}^*, A) = \Cd(F_{\fp}^*, A)/A$. Define  
\begin{equation}
\label{cycle5}
s\colon A \,\lra \, \widetilde{E}, \qquad a\mapsto \alpha(1_{\cO_{\fp}}\otimes a)-\beta(1_{\cO_{\fp}} \cdot a).
\end{equation}
A simple calculation shows that \eqref{cycle5} is $F_{\fp}^*$-equivariant (i.e.\ the image of $s$ is $F_{\fp}^*$-invariant). Also, the composition $\gamma\circ s$ is equal to the inclusion $A\to \Cd(F_{\fp}^*, A)$. Hence the pull-back of \eqref{cycle2} to $A\subseteq \Cd(F_{\fp}^*, A)$ has a canonical splitting. Consequently \eqref{cycle2} induces a short exact sequence of $F_{\fp}^*$-modules 
\begin{equation}
\label{cycle6}
\begin{CD}
0@>>> \Ccd(F_{\fp},A)@>>> E @>>> \barCd(F_{\fp}^*, A) @>>>0
\end{CD}
\end{equation} 
where $E=\widetilde{E}/s(A)$. 

\begin{remark} 
\rm Note that \eqref{cycle6} splits as a sequence of abelian groups (but not as $F_{\fp}^*$-modules). In fact \eqref{cycle0}, hence \eqref{cycle2}, and $0\to A \to \Cd(F_{\fp}^*, A)\to \barCd(F_{\fp}^*, A)\to 0$ are split exact. So if $\sigma_1: \barCd(F_{\fp}^*, A)\to \Cd(F_{\fp}^*, A)$ and $\sigma_2: \Cd(F_{\fp}^*, A)\to \widetilde{E}$ are respective splittings and if $\pi: \widetilde{E}\to E$ is the projection then $\pi\circ \sigma_2\circ \sigma_1$ is a splitting of \eqref{cycle6}.
\end{remark}

Now we return to the homomorphism $\psi\colon F_\fp^* \to A$.  Since the function $x\psi- \psi$ for $x\in F_{\fp}^*$ is equal to the constant map $y \mapsto -\psi(x)$, the class $\barpsi$ of $\psi$ in $\barCd(F_{\fp}^*, A)$ is $F_{\fp}^*$-invariant. Then applying the connecting homomorphism 
$\delta\colon H^0(F_{\fp}^*, \barCd(F_{\fp}^*, A))\to H^1(F_{\fp}^*,\Ccd(F_{\fp},A))$
induced by \eqref{cycle6}, we have
\begin{equation}
\label{cyclehom}
c_{\psi}=\delta(\barpsi)\in H^1(F_{\fp}^*, \Ccd(F_{\fp}, A)).
\end{equation} 

\subsection{Two commutative diagrams}
Let $H$ be a subgroup of $F_{\fp}^*$ and $M$ an $H$-module. Tensoring \eqref{cycle6} with $M$ yields a short exact sequence of $H$-modules
\begin{equation}
\label{cycle7}
\begin{CD}
0\lra \Cd(F_{\fp},A)\otimes M \lra E\otimes M \lra \barCd(F_{\fp}^*, A)\otimes M \lra 0
\end{CD}
\end{equation} 
(since \eqref{cycle6} splits as a sequence of abelian groups). 
Consider the following diagram:
\begin{equation} \label{e:commdiagram} \begin{gathered}
\begin{CD}
H^i(H, \Cd(F_{\fp}^*, A)\otimes M) @> x\mapsto c_{\fp} \cup x >> H^{i+1}(H, \Ccd(F_{\fp}^*, A)\otimes M)\\
@VVV @VVV\\
H^i(H, \barCd(F_{\fp}^*, A)\otimes M) @>\delta >> H^{i+1}(H, \Ccd(F_{\fp}, A)\otimes M).
\end{CD}
\end{gathered}\end{equation}
Here the cup-product is induced by the pairing \eqref{multfunctions2}, i.e.\ by 
\[ C_c(F_{\fp}^*,\bZ)\times \Cd(F_{\fp}^*, A)\to \Ccd(F_{\fp}^*, A).\] The first vertical map is induced by the projection $\Cd(F_{\fp}^*, A) \to \barCd(F_{\fp}^*, A)$, the second by $\Ccd(F_{\fp}^*, A)\to \Ccd(F_{\fp}, A), f \mapsto f_!$ and $\delta$ is the connecting homomorphism associated to \eqref{cycle7}. 
The following result follows immediately from the definitions and standard functorial properties of the cup-product.

\begin{lemma} 
\label{lemma:cupcompatibility}
\begin{enumerate}
\item[(a)] The  diagram $(\ref{e:commdiagram})$ commutes.
\item[(b)] Let $\psi: F_{\fp}^*\to A$ be a continuous homomorphism and let $\barpsi\in H^0(H, \barCd(F_{\fp}^*, A))$ be as above. We have $\delta(\barpsi \cup x) \, =\, c_{\psi} \cup x$ for all $x\in H^i(H, M)$.
\end{enumerate}
\end{lemma}

Assume now that $S$ is an arbitrary finite set of nonarchimedean places of $F$. The pairing \eqref{delta1} induces a cap-product pairing
\begin{equation*}
\cap: H^i(F^*, \cCd(S, A))\times H_j(F^*, \cC_c(\emptyset, \bZ)) \to H_{j-i}(F^*,  \cCcd(S, A))
\end{equation*}
for $j\ge i\ge 0$. In particular, taking the cap-product with $\vartheta$ yields a map 
\begin{equation}
\label{delta2}
H^0(F^*, \cCd(S, A)) \,\stackrel{\wcdot \cap \vartheta}{\lra} \, H_{n-1}(F^*, \cCcd(S, A)).
\end{equation}
It is $\bA_F^*/U^S$-equivariant since $\vartheta$ is $\bA_F^*$-invariant.

More generally let $S_1$ and $S_2$ be disjoint subsets of $S$ with $S= S_1 \cup S_2$.  Let $r$ be the cardinality of $S_1$ and choose an ordering $\fp_1, \ldots, \fp_r$ of the elements of $S_1$. Again we have a canonical cap-product pairing
\begin{equation*}
H^i(F^*, \cCd(S_2, A)^{S_1})\times H_j(F^*, \cC_c(\emptyset, \bZ)^{S_1}) \stackrel{\cap}{\lra} H_{j-i}(F^*, \cCcd(S_2, A)^{S_1})
\end{equation*}
for $j\ge i\ge 0$ so we get a map 
\begin{equation}
\label{deltaS2}
H^0(F^*, \cCd(S_2, A)^{S_1}) \,\lra \, H_{n + r-1}(F^*, \cCcd(S_2, A)^{S_1}), \, \varphi \mapsto \varphi\cap \vartheta^{S_1}.
\end{equation}
We introduce the following generalizations of \eqref{adelicfunctions}
\begin{equation}
\label{adelicfunctions2}
\cC_{?}^{\sharp}(S_1,S_2, A) \, = \, C_{?}^{\sharp}(\prod_{\fp\in S_1} F_{\fp} \times (\bA_F^{S_1})^*/ U^{S_1\cup S_2}, A)
\end{equation}
for $?\in \{\emptyset, c\}$ and $\sharp\in \{\emptyset, \diamond\}$. Also we sometimes write $\cC_{?}^0(S_1,S_2, A)$ to stress if $A$ is equipped with the discrete topology.
Again the groups \eqref{adelicfunctions2} carry a canonical $\bA_F^*$-action. We have the inclusion $\cC_c^{?}(S, A) \subseteq \cC_c^{?}(S_1, S_2, A)$ for $?\in \{\emptyset, \diamond\}$ (the map \eqref{extensionsbyzero}) and we have equality if $S_1 =\emptyset$. 

Consider the following diagram:
\begin{equation} \label{e:commdiagram2}  \begin{gathered}
\begin{CD}
H^0(F^*, \cCd(S_2, A)^{S_1}) @> \eqref{deltaS2} >> H_{n+r-1}(F^*, \cCcd(S_2, A)^{S_1})\\
@VV  \incl V @VVV\\
H^0(F^*, \cCd(S, A)) @>>>  H_{n-1}(F^*, \cCcd(S_1, S_2, A))
\end{CD}
\end{gathered}
\end{equation}
Here the first vertical map is induced by $\pi^*\colon \cCd(S_2, A)^{S_1}\to \cCd(S, A)$, where \[ \pi\colon  \bA_F^*/U^S \to (\bA_F^{S_1})^*/U^S\] denotes the projection. The second vertical map is
\begin{eqnarray*}
&& H_{n+r-1}(F^*, \cCcd(S_2, A)^{S_1}) \,\lra \,  H_{n-1}(F^*, \cCcd(S_1, S_2, A)), \\
&& \hspace{1cm} \alpha\mapsto (-1)^{r(r-1)/2}(c_{\ord_{\fp_1}} \cup \ldots \cup c_{\ord_{\fp_r}})\cap \alpha\nonumber
\end{eqnarray*}
and the lower horizontal map is the composite of \eqref{delta2} with 
\begin{equation*}
H_{n-1}(F^*, \cC_c^0(S, A)) \,\stackrel{\eqref{extensionsbyzero}_*}{\lra}\, H_{n-1}(F^*, \cC_c^0(S_1 , S_2, A)).\nonumber
\end{equation*}

The following result follows immediately from Lemma \ref{lemma:changeS} and Remark~\ref{remark:classord}.

\begin{lemma} 
\label{lemma:capord}
The  diagram $(\ref{e:commdiagram2})$ commutes.
\end{lemma}

\subsection{Homology classes associated to idele class characters} Let $R$ be a locally profinite ring and let $\chi: \bA_F^*\to R^*$ be a continuous homomorphism.  In our applications, $R$ will be a quotient of the group ring associated to the Galois group of a finite extension of $F$ and $\chi$ will be the Artin reciprocity map attached to the extension.

For a place $v$ of $F$ we denote by $\chi_v\colon F_v^* \to R^*$ the local components of $\chi$, i.e.\ the restriction to the subgroup $F_v^* \cong \{(x_{w})\in \bA_F^*\mid x_w=1 \,\,\mbox{for all $w\ne v$}\,\}$. If $\Sigma$ is a finite set of places of $F$ we also denote the restriction of $\chi$ to $(\bA_F^{\Sigma})^*= \{(x_{w})\in \bA_F^*\mid x_w=1\,\, \mbox{for all $w\in \Sigma$}\,\}$ by $\chi^{\Sigma}$. We assume that $\chi$ is unramified outside some finite set of nonarchimedean places $S$ of $F$, i.e.\ the kernel of $\chi_v$ contains $U_v$ for all $v\not\in S$. Therefore we have $\chi((x_v)_v) = \prod_v \chi_v(x_v)$ for all $(x_v)_v\in \bA_F^*$. Furthermore we assume that $\chi$ is trivial on principal ideles so that we can view $\chi$ as a homomorphism $\chi\colon \bA_F^*/F^*U^S\to R^*$ or as an element of $H^0(F^*, \cCd(S, R))$.

Assume now that $\fa$ is a closed ideal of $R$ and put $\barR = R/\fa$. Also for any $R$-module $M$ we put $\barM= \barR\otimes_R M$. We denote by $\barchi\colon \bA_F^*\to R^*\stackrel{\pr}{\lra}\barR^*$ the reduction of $\chi$ modulo $\fa$. Let $S_1, S_2$ be disjoint subsets of $S$ with $S_1\cup S_2 =S$. 
For $v\in S_1$ let $\fa_v$ be a closed ideal of $R$ contained in $\fa$ so that $\chi_v(x) \equiv 1\mod\fa_v$ for all $x\in F_v^*$ and $v\in S_1$. In particular we have $\barchi_v = 1$ for all $v\in S_1$ hence $\barchi^{S_1}$ factors through $(\bA_F^{S_1})^*/F^*U^S\to \barR^*$. Thus we can view $\barchi^{S_1}$ as an element of $H^0(F^*, \cCd(S_2, \barR)^{S_1})$.

Suppose that $S_1 = \{ v_1, \ldots, v_r \}$ and write $\fa_i=\fa_{v_i}$ for $i=1, \ldots, r$. Since the image of $\chi_{v_i}$  is contained in $1+\fa_i$ we  can consider the homomorphism 
\begin{equation*}
\der\!\chi_{v_i}\colon F_{v_i}^*\, \lra\, \barfa_i, \,\, x\mapsto  \chi_{v_i}(x)-1 \mod \fa\fa_i
\end{equation*}
and the associated cohomology class $c_{\der\!\chi_{v_i}}\in H^1(F_{v}^*, C_c(F_{v_i}, \barfa_i))$. 

Suppose now that $\fa\cdot \fa_1\cdot \ldots \cdot \fa_r=0$. Multiplication in $R$ induces a $r+1$-multilinear map 
\begin{equation*}
\mult\colon \barR\times \barfa_1\times \ldots \times \barfa_r \,\lra \, \overline{\fa_1 \ldots \fa_r} \,=\,\fa_1\ldots  \fa_r
\end{equation*}
hence a homomorphism (compare \eqref{multfunctions3})
\begin{equation*}
\bigotimes_{i=1}^r \Ccd(F_{v_i}, \barfa_i)\otimes \cCcd(S_2, \barR)^{S_1} \, \lra\, \cCcd(S_1,S_2, \fa_1\ldots  \fa_r).
\end{equation*}
This map induces a cap-product so we obtain a map
\begin{eqnarray*}
&& H_{n+r-1}(F^*, \cC_c(S_2, \barR)^{S_1}) \, \lra\, H_{n-1}(F^*, \cC_c(S_1,S_2, \fa_1\ldots  \fa_r)),\\
&& \hspace{2cm} \alpha\mapsto (c_{\der\!\chi_{v_1}}\cup\ldots \cup c_{\der\!\chi_{v_r}})\cap \alpha.
\end{eqnarray*}
Let $\iota\colon \cC_c(S_1,S_2, \fa_1\ldots  \fa_r) \hookrightarrow \cC_c(S_1,S_2, R)$ denote the inclusion.

\begin{prop} 
\label{prop:higherordervanishing}
Let $\kappa\in H_{n-1}(F^*, \cCcd(S_1,S_2, R))$ denote the image of $\chi$ under 
\begin{equation*}
H^0(F^*, \cCd(S, R)) \stackrel{\eqref{delta2}}{\lra} H_{n-1}(F^*, \cCcd(S, R))\stackrel{\eqref{extensionsbyzero}_*}{\lra} \, H_{n-1}(F^*, \cCcd(S_1 , S_2, R))
\end{equation*}
and let $\barkappa\in H_{n+r-1}(F^*, \cC_c(S_2, \barR)^{S_1})$ be the image of $\barchi^{S_1}$ under \eqref{deltaS2}. Then we have
\begin{equation*}
\kappa \, = \, \iota_*((c_{\der\!\chi_{v_1}}\cup\ldots \cup c_{\der\!\chi_{v_r}})\cap\barkappa).
\end{equation*}
In particular $\kappa=0$ if $\fa_1\ldots \fa_r=0$ holds in $R$.  
\end{prop}

{\em Proof.} The proof is based on the commutativity of a large diagram
\begin{equation}
\label{bigdiagram}
\begin{CD}
H^0(F^*, \cC_1) @> 1 >> H^r(F^*, \cC_2) @> 2 >> H_{n-1}(F^*, \cC_3)\\
@VV 3 V @VV 4 V @VV 5 V\\
H^0(F^*, \cC_4) @> 6 >> H^r(F^*, \cC_5)@> 7 >> H_{n-1}(F^*, \cC_6)\\
@AA 8 A @AA 9 A@AA 10 A\\
H^0(F^*, \cC_7) @>11 >> H^r(F^*, \cC_8)@> 12 >> H_{n-1}(F^*, \cC_9)\\
@VV 13 V @VV 14 V@A A  = A\\
H^0(F^*, \cC_{10}) @>15 >> H^r(F^*, \cC_{11})@> 16 >> H_{n-1}(F^*, \cC_{12}).
\end{CD}
\end{equation}
We are going to introduce the groups and maps involved. The $F^*$-modules in the first two columns are 
\begin{align*}
 \cC_1= \left(\bigotimes_{i=1}^r\, \Cd(F_{v_i}^*, R)\right)\otimes \cCd(S_2, R)^{S_1} & & \cC_2= \left(\bigotimes_{i=1}^r\, \Ccd(F_{v_i}^*, R)\right)\otimes \cCd(S_2, R)^{S_1}\\
 \cC_4= \left(\bigotimes_{i=1}^r\, \barCd(F_{v_i}^*, R)\right)\otimes \cCd(S_2, R)^{S_1} & & \cC_5= \left(\bigotimes_{i=1}^r\, \Ccd(F_{v_i}, R)\right)\otimes \cCd(S_2, R)^{S_1}\\
\cC_7= \left(\bigotimes_{i=1}^r\, \barCd(F_{v_i}^*, \fa_i)\right)\otimes \cCd(S_2, R)^{S_1} & &\cC_8= \left(\bigotimes_{i=1}^r\, \Ccd(F_{v_i}, \fa_i)\right)\otimes \cCd(S_2, R)^{S_1}\\
\cC_{10}= \left(\bigotimes_{i=1}^r\, \barCd(F_{v_i}^*, \barfa_i)\right)\otimes \cCd(S_2, \barR)^{S_1} & & \cC_{11}= \left(\bigotimes_{i=1}^r\, \Ccd(F_{v_i}, \barfa_i)\right)\otimes \cCd(S_2, \barR)^{S_1}
\end{align*}
(all tensor products over $R$). For the third column we put $\cC_3 = \cCd(S, R)$, $\cC_6 = 
\cCd(S_1, S_2, R)$ and $\cC_9 = \cCd(S_1, S_2, \fa_1\ldots  \fa_r) = \cC_{12}$.

The map $3$ is induced by the projections $\Cd(F_{v_i}^*, R)\to \barCd(F_{v_i}^*, R)$ and the maps $4$ and $5$ by $\Ccd(F_{v_i}^*, R)\hookrightarrow \Ccd(F_{v_i}, R), f\mapsto f_!$ for $i=1, \ldots, r$. The maps $8$, $9$ are induced by the inclusions 
$\fa_i\hookrightarrow R$, $i=1,\ldots,r$ and $10$ by $\fa_1\ldots  \fa_r\hookrightarrow R$. Similarly, the maps $13$ and $14$ are induced by the projections $\fa_i\to \barfa_i$, $R\to \barR$. The first horizontal map is given by an iterated cup-product $x\mapsto c_{v_1}\cup \ldots \cup c_{v_r} \cup x$ and the maps $6$, $11$ and $15$ are compositions of connection homomorphisms $\delta_{v_1}\circ\ldots \circ \delta_{v_r}$ as defined in Lemma \ref{lemma:cupcompatibility} (a). Finally, the maps $2$, $7$, $12$ and $16$ are all given by taking the cap-product with the $\vartheta^{S_1}\in H_{n+r-1}(F^*, \cC_c(\emptyset, \bZ)^{S_1})$. For example the paring $\cC_c(\emptyset, \bZ)^{S_1}\times \cCd(S_2, R)^{S_1}\to \cCcd(S_2, R)^{S_1}$ induces $\cC_c(\emptyset, \bZ)^{S_1}\times \cC_2 \to \cC_3$ hence 
\begin{equation*}
\cap: H^r(F^*, \cC_2) \times H_{n+r-1}(F^*, \cC_c(\emptyset, \bZ)^{S_1})\to  H_{n-1}(F^*, \cC_3).
\end{equation*}
The commutativity of the upper left square in \eqref{bigdiagram} follows from Lemma \ref{lemma:cupcompatibility} (a). All other squares are commutative for obvious reasons.

Put $\chi_i= \chi_{v_i}$, $\phi_i = \chi_i -1$ and $\psi_i = \der\!\chi_{v_i} = \phi_i \pmod{\fa\fa_i}$ for $i=1, \ldots, r$. By assumption we have $\phi_i \in \Cd(F_{v_i}, \fa_i)$. Let $\barphi_i$ and $\barpsi_i$ be the classes of $\phi_i$ and $\psi_i$ in $\barCd(F_{v_i}^*, \fa_i)$ and $\barCd(F_{v_i}^*, \barfa_i)$, respectively. We can view $\chi$ as an element of $H^0(F^*, \cC_1)$ (by identifying $\chi$ with $\bigotimes_{i=1}^r \chi_i \otimes \chi^{S_1}\in \cC_1$). By Lemma \ref{lemma:changeS} the image of $\chi$ under the composition of the maps $1$, $2$ and $5$ in \eqref{bigdiagram} is $\kappa$.

Consider the element $\eta = \bigotimes_{i=1}^r \barphi_i \otimes \chi^{S_1}\in \cC_7$. We show that $\eta$ is  $F^*$-invariant. For $x\in F^*$ we have $x \cdot \chi_i = \chi_i(x)^{-1} \chi_i$ hence $x \cdot \barphi_i = \chi_i(x)^{-1} \barphi_i$ since 
$\chi_i(x)^{-1} \barphi_i = \chi_i(x)^{-1}\chi_i - \chi_i(x)^{-1} \equiv x \cdot \chi_i - 1$ modulo constants, and the last term is $x \cdot \barphi_i $.
It follows that
\begin{equation*}
x\cdot \eta = \bigotimes_{i=1}^r (x \cdot\barphi_i) \otimes (x \cdot\chi^{S_1}) = \left(\prod_{i=1}^r \chi_i(x^{-1}) \cdot \chi^{S_1}(x^{-1})\right)\cdot \eta = \eta 
\end{equation*}
for all $x\in F^*$.  Therefore $\eta$ defines an element of $H^0(F^*, \cC_7)$, and it is clear that its image under $8$ equals the image of $\chi$ under 3.
 By the commutativity of the upper two squares in \eqref{bigdiagram} the image of $\eta$ under the composite of the maps $8$, $6$ and $7$ is
 therefore $\kappa$. 
 
 On the other hand $\eta$ is mapped under $13$ to $\bigotimes_{i=1}^r \barpsi_i \otimes \barchi^{S_1}$ so we can apply Lemma \ref{lemma:cupcompatibility} (b) to compute the image of $\eta$ under the composite of $11$, $12$ and $10$. 
Indeed, we have $\barpsi_i\in H^0(F^*, \barCd(F_{v_i}^*, \barfa_i))$ for $i=1, \ldots, r$ and $\barchi^{S_1}\in H^0(F^*,  \cCd(S_2, \barR)^{S_1})$ so by Lemma \ref{lemma:cupcompatibility} (b) $\eta$ is mapped to 
\begin{equation*}
((c_{\psi_1}\cup\ldots \cup c_{\psi_r})\cup \barchi^{S_1})\cap \vartheta^{S_1} \, =\, (c_{\psi_1}\cup\ldots \cup c_{\psi_r})\cap\partial^{S_1}(\barchi^{S_1})
\end{equation*}
under the composite of $13$, $15$ and $16$.
\enddemo

\begin{remark} 
\label{remark:cyclchar}
\rm We point out that Prop.\ \ref{prop:higherordervanishing} generalizes Prop.\ 4.6 of \cite{spiess1} (also the proof given above is simpler and more transparent than in \cite{spiess1}). To see this, let $p$ be a prime and choose a decomposition $S_p=S_1\cup S_2$ into disjoint subsets. Let $F^{\cycl}/F$ be the cyclotomic $\bZ_p$-extension of $F$. We view $\Gal(F^{\cycl}/F)$ as a subgroup of $\Gamma = 1 +2p \bZ_p\subseteq \bZ_p^*$. Let $C^{\an}(\bZ_p, \bZ_p)$ be the $\bZ_p$-algebra of locally analytic functions $\bZ_p\to \bZ_p$. An element $\gamma\in \Gamma$ defines a function $\iota(\gamma)\in C^{\an}(\bZ_p, \bZ_p)$ given by $\iota(\gamma)(s) = \gamma^s = \exp(s\log_p(\gamma))$ for $s\in \bZ_p$. Put $R = \bQ_p[X]/(X)^{r+1}$, $\fa= \fa_1=\ldots =\fa_r =(X)/(X^{r+1})$ and $\barX = X+(X^{r+1})\in R$ (where $r=\sharp(S_1)$). Define a ring homomorphism 
\begin{equation*}
\Ta_{\le r}: C^{\an}(\bZ_p, \bZ_p) \, \lra \, R,\,\,  f\mapsto \Ta_{\le r} f = \sum_{k=0}^r \frac{f^{(k)}(0)}{k!} \barX^k.
\end{equation*}
If we apply Prop.\ \ref{prop:higherordervanishing} to $S=S_p$ and the homomorphism 
\begin{equation*}
\begin{CD}
\bA_F^*/F^*U^p @> \rec >> \Gal(F^{\cycl}/F) \, \subseteq \, \Gamma @> \iota >> C^{\an}(\bZ_p, \bZ_p)^*@> \Ta_{\le r}>>R^*
\end{CD}
\end{equation*}
then it is easy to see that we essentially recover Prop.\ 4.6 of \cite{spiess1}.
\end{remark}

\subsection{Archimedean Places} Next we formulate and prove a variant of Prop.\ \ref{prop:higherordervanishing} that takes into account the infinite places of $F$. For this we fix a subset $S_3$ of $S_{\infty}$ of cardinality $s\ge 1$ and continuous homomorphisms  $\ep_v\colon F_v^*\to \{\pm 1\}=\bZ^*$ for every $v\in S_3$. We also put $F_{S_3} = \prod_{v\in S_3} F_v$ and define $\ep = \ep_{S_3}: F_{S_3}^*\to \{\pm 1\}, (x_v)_{v\in S_3}\mapsto \prod_{v\in S_3} \ep_v(x_v)$. For a subgroup $H\subseteq F_{S_3}^*$ and an $H$-module $M$ we put $M(\ep) = M\otimes \bZ(\ep)$, i.e.\ we have $M(\ep) = M$ as an abelian group and the $H$-action is given by $x\cdot m = \ep(x) x m$ for $x\in H$ and $m\in M$.

Let $S$ be a finite set of nonarchimedean places and let $S_1$ and $S_2$ be disjoint subsets of $S$ with $S= S_1 \cup S_2$. As before we write $S = \{v_1, \ldots, v_r\}$. Moreover let $A$ be a locally profinite abelian group. Since $C(F_{S_3}^*/U_{S_3}, \bZ) = \Ind^{F_{S_3}^*}_{U_{S_3}} \bZ$ we have isomorphisms of $(\bA_F^{S_1})^*$-modules
\begin{eqnarray*}
&& \cCd(S_2, A)^{S_1}\, \cong\, C(F_{S_3}^*/U_{S_3}, \bZ) \otimes \cCd(S_2, A)^{S_1\cup S_3},\\
&& \cCcd(S_2, A)^{S_1}\, \cong \, C(F_{S_3}^*/U_{S_3}, \bZ) \otimes \cCcd(S_2, A)^{S_1\cup S_3}.
\end{eqnarray*}
Thus by tensoring the ($F_{S_3}^*$-equivariant) homomorphism 
\begin{equation}
\label{augep}
C(F_{S_3}^*/U_{S_3}, \bZ) \,\lra\, \bZ(\ep), \,\, f\mapsto \sum_{x\in F_{S_3}^*/U_{S_3}} \ep(x) f(x)
\end{equation}
with $\id_{\cCd(S_2, A)^{S_1\cup S_3}}$ and $\id_{\cCcd(S_2, A)^{S_1\cup S_3}}$ respectively we obtain $(\bA_F^{S_1})^*$-equivariant maps 
\begin{eqnarray}
\label{augep2}
\cCd(S_2, A)^{S_1} \,\lra\, \cCd(S_2, A)^{S_1\cup S_3}(\ep),\\
\label{augep3}
\cCcd(S_2, A)^{S_1} \,\lra\, \cCcd(S_2, A)^{S_1\cup S_3}(\ep).
\end{eqnarray}
For $\ep =1$ the map \eqref{augep} induces a homomorphism 
\begin{equation}
\label{augep4}
H_{n-1}(F^*, \cC_c(\emptyset, \bZ)^{S_1})\lra H_{n-1}(F^*, \cC_c(\emptyset, \bZ)^{S_1\cup S_3})
\end{equation}
and we denote the image of $\vartheta^{S_1}$ under \eqref{augep4} by $\vartheta^{S_1\cup S_3}$. Again the pairing \[ \cC_c(\emptyset, \bZ)^{S_1\cup S_3}\times \cCd(S_2, A)^{S_1\cup S_3}(\ep)\to 
\cCcd(S_2, A)^{S_1\cup S_3}(\ep) \] induces a cap-product pairing so in particular a map
\begin{equation}
\label{delta2infty}
\begin{CD} 
H^0(F^*\!, \cCd(S_2, A)^{S_1\cup S_3}\!(\ep))\!\! @> \wcdot \cap \vartheta^{S_1\cup S_3} >> \!\! H_{n+r -1}(F^*\!, \cCcd(S_2, A)^{S_1\cup S_3}\!(\ep)).
\end{CD}
\end{equation}

As before we consider a locally profinite ring $R$ together with closed ideals $\fa_{v_i} =\fa_i\subseteq \fa\subseteq R$ for $i=1, \ldots, r$ and a continuous homomorphism $\chi\colon \bA_F^*/F^*U^S\to R^*$ such that $\chi_{v_i} \equiv 1\mod \fa_i$ for all $=1,\ldots, r$. Moreover assume that for every $v\in S_3$ we have given a closed ideal $\fa_v\subseteq \fa$ so that $\chi_v(-1) \equiv -\ep_v(-1)\mod \fa_v$ for all $v\in S_3$. This implies that
\begin{equation*}
\psi_v\, =\, \ep_v(1)\chi_v(1) + \ep_v(-1)\chi_v(-1) =1 + \ep_v(-1)\chi_v(-1)
\end{equation*}
lies in $\fa_v$ for all $v\in S_3$. We denote by $\barpsi_v\in \barfa_v$ the residue class of $\psi_v$ modulo $\fa\fa_v$. Here, as before, we write $\barR=R/\fa$, $\barfa_i=\fa_i/\fa\fa_i$, $\barfa_v=\fa_v/\fa\fa_v$ etc. 

Multiplication in $R$ induces a map 
\begin{equation}
\label{gradedmult4}
\bigotimes_{v\in S_3} \barfa_v \otimes_R \cCd(S_2, \barR)^{S_1\cup S_3} \,\lra\, \cCd(S_2, \overline{\prod_{v\in S_3}\fa_v})^{S_1\cup S_3}.
\end{equation}
We let $\wchi^{S_1}$ be the image of $\bigotimes_{v\in S_3} \barpsi_v\otimes \barchi^{S_1\cup S_3}$ under \eqref{gradedmult4}. Since $\chi_v(x)\psi_v = \ep_v(x)\psi_v$ we also have 
\begin{equation}
\label{psiinftyinv}
\chi_v(x)\barpsi_v = \ep_v(x)\barpsi_v
\end{equation}
for all $x\in F_v^*$ and $v\in S_3$. Therefore a simple computation shows that we have \[ \wchi^{S_1}\in H^0(F^*, \cCd(S_2, \overline{\prod_{v\in S_3}\fa_v})^{S_1\cup S_3}(\ep)). \] 

Assume that $\fa\cdot \prod_{v\in S_1\cup S_3}\fa_v= 0$. Let $\kappa\in H_{n-1}(F^*, \cCcd(S_1,S_2, R)^{S_3}(\ep))$ denote the image of $\chi$ under the composite map
\begin{eqnarray*}
H^0(F^*, \cCd(S, R)) & \stackrel{\eqref{delta2}}{\lra} & H_{n-1}(F^*, \cCcd(S, R))\\ &\stackrel{\eqref{augep3}_*}{\lra} & H_{n-1}(F^*, \cCcd(S, R)^{S_3}(\ep)) \nonumber \\
& \stackrel{\eqref{extensionsbyzero}_*}{\lra}&
H_{n-1}(F^*, \cCcd(S_1 , S_2, R)^{S_3}(\ep))\nonumber
\end{eqnarray*}
and let $\barkappa\in H_{n+r-1}(F^*, \cC_c(S_2, \overline{\prod_{v\in S_3}\fa_v})^{S_1\cup S_3}(\ep))$ be the image of $\wchi^{S_1}$ under \eqref{delta2infty}. Let \[ \iota: \cC_c(S_1,S_2, \prod_{v\in S_1\cup S_3}\fa_v)^{S_3} \hookrightarrow \cC_c(S_1,S_2, R)^{S_3} \] denote the inclusion.

\begin{prop} 
\label{prop:higherordervanishing2}
We have
\begin{equation}
\label{higherordervanishing4}
\kappa \, = \, \iota_*((c_{\der\!\chi_{v_1}}\cup\ldots \cup c_{\der\!\chi_{v_r}})\cap\barkappa).
\end{equation}
In particular $\kappa=0$ if $\prod_{v\in S_1\cup S_3}\fa_v=0$ holds in $R$. 
\end{prop}

{\em Proof.} Let $\wchi$ be the image of $\bigotimes_{v\in S_3} \psi_v\otimes  \chi^{S_3}$ under the canonical map
\begin{equation*}
\bigotimes_{v\in S_3} \fa_v \otimes_R \cCd(S, R)^{S_3}(\ep) \,\lra\, \cCd(S, \prod_{v\in S_3}\fa_v)^{S_3}(\ep).
\end{equation*}
Again using \eqref{psiinftyinv} we obtain $\wchi\in H^0(F^*, \cCd(S, \prod_{v\in S_3}\fa_v)^{S_3}(\ep))$. We have a commutative diagram 
\begin{equation}
\label{diagramkey}
\begin{CD}
H^0(F^*, \cCd(S, R)) @>\eqref{delta2}>> H_{n-1}(F^*, \cCcd(S, R))\\
@VV \eqref{augep2}_* V @VV\eqref{augep3}_* V\\
H^0(F^*, \cCd(S, R)^{S_3}(\ep)) @ > \eqref{delta2infty} >> H_{n-1}(F^*, \cCcd(S, R)^{S_3}(\ep))\\
@AA \incl A@AAA\\
H^0(F^*, \cCd(S, \prod_{v\in S_3}\fa_v)^{S_3}(\ep)) @ > \eqref{delta2infty}  >> H_{n-1}(F^*, \cCcd(S, \prod_{v\in S_3}\fa_v)^{S_3}(\ep))
\end{CD}
\end{equation}
where the lower vertical maps are induced by the inclusion $\prod_{v\in S_3}\fa_v\subseteq R$. The image of $\chi$ under $\eqref{augep2}_*$ is equal to $\wchi\in H^0(F^*, \cC(S, \prod_{v\in S_3}\fa_v)^{S_3}(\ep))$. Thus it suffices to show that the image of $\wchi$ under the composition of the lower horizontal map in \eqref{diagramkey} with 
\[ H_{n-1}(F^*, \cCcd(S, \prod_{v\in S_3}\fa_v)^{S_3}(\ep))\stackrel{\eqref{extensionsbyzero}_*}{\lra}H_{n-1}(F^*, \cCcd(S_1 , S_2, \prod_{v\in S_3}\fa_v)^{S_3}(\ep)) \] agrees with the right hand side of \eqref{higherordervanishing4}. For that we can follow the proof of Prop.\ \eqref{prop:higherordervanishing2}. We just have to modify diagram \eqref{bigdiagram} accordingly. For example $\cC_1$ should be replaced by 
\[ \cC_1'= \bigotimes_{i=1}^r\Cd(F_{v_i}^*, R)\otimes \cCcd(S, \prod_{v\in S_3}\fa_v)^{S_1\cup S_3}\!(\ep),\] $\cC_2$ by \[ \cC_2'= \left(\bigotimes_{i=1}^r\Ccd(F_{v_i}^*, R)\right)\otimes \cCcd(S, \prod_{v\in S_3}\fa_v)^{S_1\cup S_3}\!(\ep),\] etc. The details will be left to the reader.\enddemo

\section{Integrality properties of the Eisenstein cocycle}
\label{section:eisenstein}

\subsection{Shintani cocycle} \label{s:shintanicocycle}

We recall briefly the construction of the {\it Shintani cocycle} as defined in \cite{cdg}. Let $V$ be an $n$-dimensional $\bQ$-vector space. Given linearly independent vectors $v_1, \ldots, v_m\in V$ we define the rational open cone $C(v_1, \ldots, v_m)\subseteq V_{\infty} = V\otimes_\bQ \bR$ as
\[
C(v_1,  \ldots, v_m)\,\, =\,\, \left\{\sum_{i=1}^m \, t_i x_i\mid \, t_i\in \bR_+  \,\,\forall \, i=1, \ldots ,m\right\}.
\]
Let $\sV$ be the $\GL(V)$-stable subset of $V_{\infty}$ consisting of vectors $w$ that are not contained in any rational hyperplane of $V_{\infty}$ (i.e\ $w\in \sV$ iff $v^*(w)\ne 0$ for any nonzero linear form $v^*: V\to \bQ$).

Let $\cK=\cK_V$ denote the subgroup of the space of functions $f\colon V_{\infty}\to \bZ$ generated by the characteristic functions $1_C$ of rational open cones. The group $\GL(V)$ acts on $\cK$ via $(g \cdot f)(v) = f(g^{-1} v)$ for $g\in \GL(V)$ and $v\in V_{\infty}$. We also consider subsets of $V_{\infty}$ of the form
\begin{eqnarray*}
L & =& \left\{\sum_{i=1}^m \, t_i x_i\mid \, t_1\in \bR, t_i\in \bR_+  \,\,\forall \, i=2, \ldots ,m\right\}\\
& =& C(v_1,  v_2, \ldots, v_m)\cup C(-v_1,  v_2, \ldots, v_m)\cup C(v_2, \ldots, v_m)
\end{eqnarray*}
(called {\it wedges}) where $v_1, \ldots, v_m\in V$ are again linearly independent. Note that $1_L(v+v_1) = 1_L(v)$ for all $v\in V$. The $\GL(V)$-stable subgroup of $\cK$ generated by characteristic functions of wedges will be denoted by $\cL$. 

Given $v_1, \ldots, v_n\in V$ and $w\in \sV$ the function $c_w(v_1, \ldots, v_n)\colon V_{\infty} \to \bZ$ is defined as $c_w(v_1, \ldots, v_n)\equiv 0$ if $v_1, \ldots, v_n$ are linearly dependent and as
\begin{equation*}
c_w(v_1, \ldots, v_n)(v) \,\,\, = \,\,\, \underset{\ep\to 0+}\lim \, 1_{C(v_1,  \ldots, v_n)}(v+ \ep w)
\end{equation*}
otherwise. In the latter case $c_w(v_1, \ldots, v_n)$ is the characteristic function of the union of $C(v_1,  v_2, \ldots, v_n)$ with some of its boundary cones (more precisely if we write $w= \sum_{i=1}^n t_iv_i$ and if $\Sigma_0\subseteq \{v_1,  v_2, \ldots, v_n\}$ is the set of $v_i$ with $t_i<0$ then $c_w(v_1, \ldots, v_n)$ is the characteristic functions of the union of all cones generated by subsets $\Sigma$ of $\{v_1,  v_2, \ldots, v_n\}$ containing $\Sigma_0$). Note that $c_{gw}(gv_1, \ldots, gv_n)(gv) = c_w(v_1, \ldots, v_n)(v)$ for all $g\in \GL(V)$. Also $c_w(v_1, \ldots, v_n) = c_{tw}(v_1, \ldots, v_n)$ for any $t\in \bR_+$ so $c_w(v_1, \ldots, v_n)$ depends only on the ray $\bR_+ w$. We denote the set of rays generated by elements of $\sV$  by $\sR$ and define $c_{Q}(v_1, \ldots, v_n) =  c_w(v_1, \ldots, v_n)$ if $Q\in \sR$ and $w\in Q$. The $\GL(V)$-operation on $\sV$ induces a $\GL(V)$-operation on $\sR$. 
 
Let $\cN= \cK/\cL$ and let $\Hom(\sR, \cN)$ be the abelian group of maps $f\colon \sR\to \cN$ endowed with the $\GL(V)$-action $(gf)(Q) = gf(g^{-1}Q)$. We also need to introduce the sign character $\sgn\colon \bR^*\to \{\pm 1\} = \bZ^*, x\mapsto \sgn(x)$. By abuse of notation we denote the homomorphism $\sgn\circ \det$ by $\sgn$ as well
\begin{equation}
\label{sign}
\sgn= \sgn\circ \det: \GL(V_{\bR}) \to \bZ^*.
\end{equation}
Similarly for a subgroup $\Gamma$ of $\GL(V_{\bR})$, we denote the restriction of \eqref{sign} to $\Gamma$ simply by $\sgn$. Fix a nonzero vector $v$ of $V$ and a determinant form $\omega$ on $V$. The {\it Shintani cocycle} is the homogeneous cocycle 
\begin{eqnarray*}
&&\Xi_{\Sh, v}: \GL(V)^n\,\lra \, \Hom(\sR,\cN)(\sgn)\\
&&\Xi_{\Sh,v}(g_1,\ldots, g_n)(Q)\,\, =\,\, \sgn(\omega(g_1v, \ldots, g_nv))\,c_{Q}(g_1v, \ldots, g_nv)\nonumber
\end{eqnarray*}
(for the proof of the cocycle relation see \cite{cdg}, Thm.\ 1.1 and 1.6). Since \[ \Xi_{\Sh,gv}(g_1,\ldots, g_n) = g\Xi_{\Sh,v}(g^{-1}g_1g,\ldots, g^{-1}g_ng) \] we see that the cohomology class of $\Xi_{\Sh,v}$ is independent of the choice of $v\in V$. We will denote this class simply by $\Xi_{\Sh}\in H^{n-1}(\GL(V), \Hom(\sR,\cN)(\sgn))$. This construction in functorial in the following sense: if $V'$ is another $n$-dimensional $\bQ$-vector space and if $\cK'$, $\cL'$, $\cN'$, $\Xi_{\Sh}'$ etc.\ denote the same objects as above for $V'$ then an isomorphism $\psi\colon V\to V'$ induces an isomorphism $\psi^*$
between the cohomology groups and we have $\psi^*(\Xi_{\Sh}')=\Xi_{\Sh}$.

\subsection{Solomon-Hu pairing and Eisenstein cocycle}

Let $\bZ\{V\}$ be the abelian group of maps $V\to \bZ$. We have $\bZ[V]\subseteq \bZ\{V\}$ and the multiplication of $\bZ[V]$ extends to a $\bZ[V]$-module structure on $\bZ\{V\}$. Consider the pairing
\begin{equation}
\label{solhu1}
\llangle\,\,\,, \,\,\, \rrangle\,:\,\, \cK \times C_c(V_\bhatZ, \bZ) \,\lra \, \bZ\{V\}
\end{equation}
given by $\llangle f, \Phi \rrangle(v) = f(v)\Phi(v)$. Let $S$ denote the multiplicative subset of $\bZ[V]$ generated by the set $\{[v]-[0]\mid v\in V, v\ne 0\}$. The image of \eqref{solhu1} composed with the localization map $\bZ\{V\}\to S^{-1}\bZ\{V\}$ is contained in $S^{-1}\bZ[V]\subseteq S^{-1}\bZ\{V\}$. In fact if $C\subseteq V_{\infty}$ is a rational cone and $\Phi\in C_c(V_\bhatZ, \bZ)$ we may choose linearly independent vectors $v_1, \ldots, v_m\in V$ with $C= C(v_1, \ldots, v_m)$ so that $\Phi$ is periodic with respect to $v_1, \ldots, v_m$, i.e.\ $\Phi(x+v_i) = \Phi(x)$ for all $i=1,\ldots,m$ and $x\in V_\bhatZ$. Then we have
\begin{equation*}
\llangle 1_C, \Phi \rrangle\, = \, \prod_{i=1}^m ([0]- [v_i])^{-1}\,\sum_{v\in V\cap P(v_1, \ldots, v_m)} \Phi(v) [v].
\end{equation*}
where $P(v_1, \ldots, v_m)$ is the ``half-open" parallelpiped
\[ P(v_1, \ldots, v_m) \,\, =\,\, \left\{\sum_{i=1}^m \, t_i x_i\mid \, 0 < t_i \le 1,  \,\,\forall \, i=1, \ldots ,m\right\}. \]

Also if $L$ is a wedge we can choose linearly independent $v_1, \ldots, v_m\in V$ with \[ L=C(v_1,  v_2, \ldots, v_m)\cup C(-v_1,  v_2, \ldots, v_m)\cup C(v_2, \ldots, v_m)\] and so that $\Phi$ is periodic with resp.\ to  $v_1, \ldots, v_m$. Then $[v_1] \cdot\llangle 1_L, \Phi \rrangle = \llangle 1_L, \Phi \rrangle$ since both $1_L$ and $\Phi$ are periodic with respect to $v_1$. Therefore \eqref{solhu1} induces a pairing
\begin{equation}
\label{solhu2}
\llangle\,\,\,, \,\,\, \rrangle\,:\,\, \cN \times C_c(V_\bhatZ, \bZ) \,\lra \, S^{-1}\bZ[V].
\end{equation}
By composing \eqref{solhu2} with the canonical map $S^{-1}\bZ[V]\to S^{-1}\bZ[\![V]\!]$ we obtain a pairing 
\begin{equation}
\label{solhu3}
\llangle\,\,\,, \,\,\, \rrangle\,:\,\, \cN \times C_c(V_\bhatZ, \bZ) \,\lra \, S^{-1}\bZ[\![V]\!].
\end{equation}
We denote by $\cR_V$ the algebra $\prod_{n\ge 0} \Sym^n V$ and let $\cQ_V$ be its quotient field. A choice of a basis $(e_1, \ldots  e_n)$ of $V$ induces isomorphisms between $\cR_V$ resp.\ $\cQ_V$ and the power series ring $\bQ[\![ z_1,\ldots, z_n]\!]$ resp.\ the field of Laurent series $\bQ(\!( z_1, \ldots, z_m)\!)$. The group homomorphism 
 $\exp\colon V\to \cR_V^*, v\mapsto \exp(v) = \sum_{i=0}^{\infty} v^n/n!$ induces  ring homomorphisms
\begin{equation}
\label{exp}
\exp\colon \bZ[\![V]\!]\lra \cR_V, \qquad \exp\colon S^{-1}\bZ[\![V]\!]\lra \cQ_V
\end{equation}
and by composing \eqref{solhu3} with \eqref{exp} we obtain the Solomon-Hu pairing 
\begin{equation}
\label{solhu4}
\llangle\,\,\,, \,\,\, \rrangle\,:\,\, \cN \times C_c(V_\bhatZ, \bZ) \,\lra \, \cQ_V.
\end{equation}
The $\GL(V)$-action on $V$ induces a $\GL(V)$-action on $\bZ[V]$, $\bZ\{V\}$, $S^{-1}\bZ[V]$, $\cR_V$ and $\cQ_V$ and 
the pairings \eqref{solhu1}, \eqref{solhu3} and \eqref{solhu4} are $\GL(V)$-equivariant.

Given abelian groups $A$ and $B$ we denote by $\sM(A, B)$ the abelian group of maps 
\begin{equation*}
\beta\colon \sR\times A \,\lra \, B,\qquad (Q,a) \mapsto \beta(Q, a)
\end{equation*}
that are linear in the second component $a\in A$. If $\Gamma\subseteq \GL(V)$ is a subgroup and $A$ and $B$ are $\Gamma$-modules then we define a $\Gamma$-action on $\sM(A,B)$ by $(\gamma\beta)(Q, a)= \gamma (\beta(\gamma^{-1}Q, \gamma^{-1} a))$.

The pairing \eqref{solhu4} induces a homomorphism of $\GL(V)$-modules 
\begin{equation}
\label{solhumod}
\Hom(\sR, \cN)\to  \sM(C_c(V_\bhatZ,\bZ), \cQ_V), \qquad f\mapsto\left( (Q,\Phi)\mapsto \llangle f(Q), \Phi \rrangle\right).
\end{equation}
The image of the Shintani cocycle $\Xi_{\Sh}$ under the map
\begin{equation}
\label{eisen}
H^{n-1}(\GL(V), \Hom(\sR, \cN)(\sgn)) \to H^{n-1}(\GL(V),\sM(C_c(V_\bhatZ,\bZ), \cQ_V)(\sgn))
\end{equation}
induced by \eqref{solhumod} is what is usually called the {\it Eisenstein cocycle} (see e.g.\ \cite{sczech, nori, daschar, cdg}). In the following we introduce an ``integral" variant of it. 

\subsection{Integrality properties}
\label{subsection:cassou}

We get rid of the denominators in $S$ by {\it smoothing} (also called {\it Cassou-Nogu{\`e}s trick}).  We fix a triple $(\ell, L, L')$ consisting of a prime number $\ell$ and $\bZ_{\ell}$-sublattices $L'\subseteq L$ of $V_{\ell}$ with $[L:L']=\ell$.  Let $\Gamma^{\sL}\subseteq \GL(V)$ be the subgroup of $\gamma\in \GL(V)$ with $\gamma L=L$ and $\gamma L'=L'$. We consider the subgroup $\cK^{\sL}$ of $\cK$ generated by characteristic functions of cones $C$ of the form $C(v_1,  \ldots, v_m)$ with $v_i\in V\cap L$ and $v_i\not\in L'$ for $i=1,\ldots, m$ and define $\cL^{\sL}$ accordingly. These are $\Gamma^{\sL}$-stable subgroups of $\cK$, so $\cN^{\sL} = \cK^{\sL}/\cL^{\sL}$ is a $\Gamma^{\sL}$-module as well. For $v\in (V\cap L)-L'$ the Shintani cocycle $\Xi_{\Sh, v}$---when restricted to $\Gamma^{\sL}$----takes values in $\cN^{\sL}$ so we obtain a homogeneous cocycle \[ \Xi_{\Sh, \sL, v}\in Z^{n-1}(\Gamma^{\sL}, \Hom(\sR,\cN^{\sL})(\sgn)) \] and a cohomology class \[ \Xi_{\Sh,\sL}\in H^{n-1}(\Gamma^{\sL}, \Hom(\sR,\cN^{\sL})(\sgn))\] as before.

We define $\phi_{\sL} \in C_c(V_{\ell}, \bZ)$ by $\phi_{\sL} = 1_L -  \ell 1_{L'}$, i.e.\
\begin{equation}
\label{cassousmooth}
\phi_{\sL}(v)  \,\,\, = \,\,\, \left\{\begin{array}{ll} 1 &  \mbox{if $v\in L-L'$,}\\
1-\ell &  \mbox{if $v\in L'$,}\\
0 & \mbox{if $v\in V_{\ell}-L$.}
\end{array}\right.
\end{equation}
Note that in $\bZ[\mu_{\ell}]$ we have 
\begin{equation}
\label{fourier}
\phi_{\sL}(v)  \,\,\, = \,\,\,  - \sum_{\psi} \psi(v)\qquad \mbox{for all $v\in L$}
\end{equation}
where the sum extends over the $\ell-1$ characters $\psi\colon L\to \mu_{\ell} \subseteq \bZ[\mu_{\ell}]^*$ with $\Ker(\psi)= L'$ (here $\mu_{\ell}\subseteq \barQ$ denotes the group of $\ell$-th roots of unity).

Since $\phi_{\sL}$ is
fixed under the action of $\Gamma^{\sL}$ the map
\begin{equation}
\label{cassou}
C_c(V_{\bhatZ^{\ell}}, \bZ) \lra C_c(V_{\bhatZ^{\ell}}, \bZ) \otimes C_c(V_{\ell},\bZ)\cong C_c(V_{\bhatZ}, \bZ),  \qquad \phi\mapsto \phi\otimes \phi_{\sL}
\end{equation}
is a homomorphism of $\Gamma^{\sL}$-modules. Hence we obtain a $\Gamma^{\sL}$-equivariant pairing
\begin{equation}
\label{solhu5}
\llangle\,\,\,, \,\,\, \rrangle^{\sL}\,\colon \,\, \cN^{\sL} \times C_c(V_{\bhatZ^{\ell}}, \bZ) \lra  \cN \times C_c(V_{\bhatZ}, \bZ) \stackrel{\eqref{solhu3}}{\lra} S^{-1}\bZ[V]
\end{equation}
where the first map is induced by the inclusion $\cK^{\sL}\hookrightarrow \cK$ in the first component and by \eqref{cassou} in the second. 

In order to get rid of denominators in \eqref{solhu4} we first look at the pairing on the level of lattices i.e.\ we consider a $\bZ$-lattice $\Lambda$ in $V$ with $\Lambda_{\ell} = L$. When we restrict the second component in \eqref{solhu5} to $C_c(\Lambda_{\bhatZ^{\ell}}, \bZ)\subseteq C_c(V_{\bhatZ^{\ell}}, \bZ)$ we obtain a $\Gamma^{\sL}_{\Lambda} = \Gamma^{\sL}\cap \GL(\Lambda)$-equivariant pairing
\begin{equation}
\label{solhu4lat}
\llangle\,\,\,, \,\,\, \rrangle\,\colon\,\, \cN^{\sL} \times C_c(\Lambda_{\bhatZ^{\ell}}, \bZ) \lra S^{-1}\bZ[\Lambda],
\end{equation}
where now $S$ denotes the multiplicative subset of $\bZ[\Lambda]$ generated by the set \[ \{[\lambda]-[0]\mid v\in \Lambda, \lambda\ne 0\}. \] After embedding $S^{-1}\bZ[\Lambda]$ into the quotient field of $\bQ[\![\Lambda]\!]$ we will show that the image of \eqref{solhu4lat} is contained in $\bZ[1/\ell][\![\Lambda]\!]$.

\begin{prop}
\label{prop:eisenint}
Let $\phi \in C_c(\Lambda_{\bhatZ^{\ell}}, \bZ)$ and let $C = C(v_1, \ldots, v_m)$ be a cone with $v_i\in V\cap L$ and $v_i\not\in L'$ for $i=1,\ldots, m$. Then we have
\[
\llangle 1_C, \phi\otimes \phi_{\sL} \rrangle\in\bZ[1/\ell][\![\Lambda]\!]. 
\]
Moreover if $\ell \ge m+2$ the constant term of $\llangle 1_C,\phi\otimes \phi_{\sL} \rrangle$ lies in $\bZ$ (i.e.\ the image under $\aug\colon\bZ[1/\ell][\![\Lambda]\!] \to \bZ[1/\ell]$).
\end{prop}

{\em Proof.} Put $\Phi = \phi\otimes \phi_{\sL}$ and $\Phi_0 = \phi\otimes 1_L\in C_c(\Lambda_{\bhatZ}, \bZ)$. By \eqref{fourier} we get in $C_c(\Lambda_{\bhatZ}, \bZ[\zeta_{\ell}])$
\begin{equation}
\label{fourier2}
\Phi(v) \,\,\, =\,\,\, -\sum_{\psi} \psi(v) \Phi_0(v)
\end{equation}  
for all $v\in \Lambda$. There exists a sublattice $\Lambda_1 \subseteq \Lambda$ such that $\Phi_0$ is periodic with respect to $\Lambda_1$. Since the local component of $\Phi_0$ at $\ell$ is $1_{L}$ we may assume that $(\Lambda_1)_{\ell}= \Lambda_{\ell} = L$ hence $\ell\nmid [\Lambda:\Lambda_1]$. By multiplying $v_1,  \ldots, v_m$ with some integer $k$ prime to $\ell$ if necessary we can assume also that $v_1,  \ldots, v_m\in \Lambda_1$ (hence $\Phi$ is periodic with respect to $\ell \sum_{i=1}^m \bZ v_i$). Since $v_i\not\in L'$ we note that $\psi(v_i)$ is a primitive $\ell$th root of unity for all characters $\psi\colon L\to\mu_{\ell}$ as above and $i=1,\ldots, m$.

Note that the parallelepiped $\wP= P(\ell v_1, \ldots, \ell v_m)$ is the disjoint union of sets of the form $(\sum_{i=1}^m  n_i v_i)+P$ with $P= P(v_1, \ldots, v_m)$ and $0\le n_1,\ldots, n_m\le \ell-1$. Because $\Phi_0$ is periodic with respect to $\sum_{i=1}^m  \bZ v_i$ by \eqref{fourier2} we have 
\begin{eqnarray*}
&& \sum_{v\in V\cap \wP} \Phi(v) [v]\,\,= \,\, - \sum_{\psi}\sum_{v\in \Lambda\cap \wP}  \psi(v) \Phi_0(v) [v]\\ 
\nonumber && = \,\,  - \sum_{\psi}\sum_{v\in \Lambda\cap P} \psi(v)\Phi_0(v)[v] \sum_{n_1,\ldots, n_m=0}^{\ell-1} \, \prod_{i=1}^m  \psi(v_i)^{n_i} [v_i]^{n_i}    \\
\nonumber && = \,\,  - \prod_{i=1}^m (1- [\ell v_i])\wcdot \sum_{\psi} \frac{\sum_{v\in \Lambda\cap P}  \psi(v)\Phi_0(v) [v]  
}{\prod_{i=1}^m(1- \psi(v_i)[v_i])}.
\end{eqnarray*}
This equality holds in the quotient field of $\bZ[\mu_{\ell}][\![\Lambda]\!]$. Note however that for $v\in \Lambda$ we have
\begin{equation*}
1 - \zeta_{\ell} [v] = (1-\zeta_{\ell}) (1 + \frac{\zeta_{\ell}}{1-\zeta_{\ell}}(1-[v])) = (1-\zeta_{\ell}) (1 + x)
\end{equation*}  
with $x= \frac{\zeta_{\ell}}{1-\zeta_{\ell}}(1-[v])$. Since $x$ lies in the augmentation ideal of $\bZ[1/\ell, \mu_{\ell}][\Lambda]$ the element $1+x$ and therefore also $1 - \zeta_{\ell} [v]$ is invertible in $\bZ[1/\ell, \zeta_{\ell}][\![\Lambda]\!]$. It follows that
\begin{equation}
\label{eisenint2}
\llangle 1_C, \Phi\rrangle\,=  \prod_{i=1}^m (1- [\ell v_i])^{-1} \sum_{v\in V\cap \wP} \Phi(v) [v] \,= -\sum_{\psi} \frac{\sum_{v\in \Lambda\cap P}  \psi(v)\Phi_0(v) [v]
}{\prod_{i=1}^m(1- \psi(v_i)[v_i])}
\end{equation}  
lies in  $\bZ[1/\ell, \zeta_{\ell}][\![\Lambda]\!]$. On the other hand $\llangle 1_C, \Phi\rrangle$ lies in the quotient field of $\bZ[1/\ell][\![\Lambda]\!]$, hence it lies in $\bZ[1/\ell][\![\Lambda]\!]$.

For the second assertion note that $\aug(\llangle 1_C, \Phi\rrangle)\in (1-\zeta_{\ell})^{-m}\bZ[\zeta_{\ell}]$ by the above computation and $(1-\zeta_{\ell})^{-m}\bZ[\zeta_{\ell}]\cap \bQ=\bZ$ if $m \le \ell-2$.\enddemo

\subsection{Puiseux series valued pairing}
According to Prop.\ \ref{prop:eisenint} the pairing \eqref{solhu4lat} takes the form 
\begin{equation}
\label{solhu4int2}
\llangle\,\,\,, \,\,\, \rrangle\,:\,\, \cN^{\sL} \times C_c(\Lambda_{\bhatZ^{\ell}}, \bZ) \, \lra\, \bZ[1/\ell][\![\Lambda]\!].
\end{equation}
In order to pass from the level of lattices to $V$ we introduce a certain type of Puiseux series ring. Let $R$ be any ring. For a lattice $\Lambda$ in $V$ the inclusion $\Lambda\hookrightarrow V$ induces an embedding $R[\Lambda]\hookrightarrow R[V]$ and $R[\![\Lambda]\!]\hookrightarrow R[\![V]\!]$ (to see that $R[\![\Lambda]\!]\to R[\![V]\!]$ is injective it is enough to see that $I(V)^n\cap R[\Lambda] = I(\Lambda)^n$ for all $n$ or that $I(\Lambda')^n\cap R[\Lambda] = I(\Lambda)^n$ for any lattice $\Lambda'$ containing $\Lambda$; the latter can be seen easily by induction on $n$). We define the subring $R\{\!\{V\}\!\}^{\sL}$ of $R[\![V]\!]$ as the union of all $R[\![\Lambda]\!]$ where $\Lambda$ ranges over lattices such that $\Lambda_\ell = L$. The group $\GL(V)$ acts on $R[\![V]\!]$ as a group of ring automorphisms and $R\{\!\{V\}\!\}^{\sL}$ is a $\Gamma^{\sL}$-stable subring. The pairing \eqref{solhu4int2} is compatible with monomorphisms
that are induced by inclusions of lattices. Hence the image of \eqref{solhu5} is contained in $\bZ[1/\ell]\{\!\{V\}\!\}^{\sL}$.

We can get rid of denominators altogether by applying smoothing twice. More precisely let $\sL = \{\sL_i\}_{i\in I}$ be a collection of triples $\sL_i= (\ell_i, L_i, L_i')$ as above. We assume that the index set $I$ is finite and has at least two elements and that $\ell_i\ne \ell_{i'}$ for $i\ne i'$. Put $T = \{\ell_i\mid i\in I\}$ and define $\Gamma^{\sL}$, $\cK^{\sL}$, $\cL^{\sL}$ as the intersection of $\Gamma^{\sL_i}$, $\cK^{\sL_i}$, $\cL^{\sL_i}$ respectively. For defining $R\{\!\{V\}\!\}^{\sL}$ we now take lattices $\Lambda$ with $\Lambda_{\ell_i} = L_i$ for all $i\in I$. Such lattices will be called compatible with the data $\sL$.
Using the $\Gamma^{\sL}$-equivariant map
\begin{equation*}
C_c(V_{\bhatZ^T}, \bZ) \lra C_c(V_{\bhatZ}, \bZ),  \,\, \phi\mapsto \phi\otimes \bigotimes_{i\in I} \phi_{\sL_i}
\end{equation*}
instead of \eqref{cassou} we obtain as before a $\Gamma^{\sL}$-equivariant pairing
\begin{equation}
\label{solhu6}
\llangle\,\,\,, \,\,\, \rrangle\,:\,\, \cN^{\sL} \times C_c(V_{\bhatZ^T}, \bZ)\, \lra \, S^{-1}\bZ[V].
\end{equation}
We obtain

\begin{corollary}
\label{corollary:lat0}
The image of \eqref{solhu6} lies in $\bZ\{\!\{V\}\!\}^{\sL}$.
\end{corollary}

Summarizing, if $\sL$ consists of just one triple $\sL=(\ell, L, L')$ the pairing \eqref{solhu5} takes the form 
\begin{equation*}
\llangle\,\,\,, \,\,\, \rrangle^{\sL}\,:\,\, \cN^{\sL} \times C_c(V_{\bhatZ^{\ell}}, \bZ) \, \lra\, \bZ[1/\ell]\{\!\{V\}\!\}^{\sL}.
\end{equation*}
If $\sL$ involves more than one prime the pairing is of the form
\begin{equation}
\label{solhu5a}
\llangle\,\,\,, \,\,\, \rrangle\,:\,\, \cN^{\sL} \times C_c(V_{\bhatZ^T}, \bZ) \, \lra\, \bZ\{\!\{V\}\!\}^{\sL}.
\end{equation}
By passing to the constant term (i.e.\ by composing it with the augmentation map) we get a pairing 
\begin{equation}
\label{solhu4b}
\llangle\,\,\,, \,\,\, \rrangle\,:\,\, \cN^{\sL} \times C_c(V_{\bhatZ^T}, \bZ)\, \lra \, \bZ.
\end{equation}
This holds even if $\sL$ consists of a single triple $(\ell, L, L')$ provided that $\ell\ge n+2$.
 
 \subsection{Locally polynomial functions}
We now introduce  another variant of  the Solomon--Hu pairing. We denote by $\Pol_V$ the ring of polynomials on $V$. More precisely $\Pol_V$ is the
symmetric algebra $\Sym^{\bu} V^{\vee}= \bigoplus_{n\in \bN} \Sym^n V^{\vee}$ of the dual $V^{\vee}$ of $V$. Any $P\in\Pol_V$ induces a polynomial function $V\to \bQ, v\mapsto P(v)$. In fact a choice of a basis $(v_1,\ldots, v_n)$ induces an isomorphism $\Pol_V\cong \bQ[X_1, \ldots, X_n], P\mapsto \widetilde{P}$ so that $P(\sum_{i=1}^n x_i v_i) = \widetilde{P}(x_1, \ldots, x_n)$. We let $\GL(V)$ act on $V^{\vee}$ via $(g\cdot \lambda)(v) = \lambda(g^{-1}v)$ for all $\lambda \in V^{\vee}$, $v\in V$ and $g\in \GL(V)$. This induces a $\GL(V)$-action on $\Pol_V$ by $(g\cdot P)(v) = P(g^{-1}v)$. 

Given $P\in \Pol_V$ the map $V\to \bQ, v\mapsto P(v)$ extends to a homomorphism of abelian groups 
\begin{equation}
\label{pol1}
P\colon \bQ[V]\, \lra\, \bQ, \, P\left(\sum_{i=1}^r a_i [v_i]\right) \, =\,  \sum_{i=1}^r a_i P(v_i)
\end{equation}

\begin{lemma} 
\label{lemma:poleval}
(a) For $P\in \Pol_V$ the homomorphism \eqref{pol1} extends uniquely to a continuous homomorphism of abelian groups
\begin{equation}
\label{pol2}
\bQ[\![V]\!]\, \lra \, \bQ, \,\,\,x\mapsto P(x).
\end{equation}

\noi (b) The pairing $\Pol_V\times \bQ[\![V]\!]\to \bQ, (P, x) \mapsto P(x)$ is $\GL(V)$-equivariant.
\end{lemma}
 
{\em Proof.} We show that $I(V)^m$ lies in the kernel of \eqref{pol1} for $m\ge \deg(P)+1$. Define 
\[
D_P: \bQ[V]\lra \bQ[V], \,  \sum_{i=1}^r a_i [v_i]\mapsto D_P\left(\sum_{i=1}^r a_i [v_i]\right) \,=\, \sum_{i=1}^r a_i P(v_i) [v_i]
\]
Then $\aug \circ D_P$ equals \eqref{pol1} and we have $D_{P_1} \circ D_{P_2} = D_{P_1 P_2}$. Also if $P= \lambda\in V^{\vee}$ is a linear form then $D_{\lambda}$ is a derivation so that $D_{\lambda}(I(V)^{m+1}) \subseteq I(V)^m$. Combining these facts we see that $D_P$ maps $I(V)^{\deg P+1}$ into $I(V)$ so \eqref{pol1} maps $I(V)^m$ to $0$ for $m\ge \deg P +1$. This proves (a). The second statement is obvious. 
\enddemo

We define the ring $\cP_c(V, \bQ)$ of locally polynomial functions on $V_{\bhatZ}$ by \[ \cP_c(V, \bQ) = C_c(V_{\bhatZ}, \bZ)\otimes \Pol_V. \] Every element $h \in \cP_c(V, \bQ)$ induces a function $V\to \bQ, v\mapsto h(v)$. Indeed if $h$ is of the form $h= \sum_i^m\Phi_i\otimes P_i\in C_c(V_{\bhatZ}, \bZ)\otimes \Pol_V$ then $v\mapsto h(v)$ is given by $h(v) = \sum_i^m\Phi_i(v)\cdot P_i(v)$. For a lattice $\Lambda\subseteq V$ we let $\Pol_V(\Lambda)$ be the subring of $\Pol_V$ consisting of polynomials $P\in \Pol_V$ with $P(v)\in \bZ$ for all $v\in \Lambda$. We can view $C_c(\Lambda_{\bhatZ}, \bZ)\otimes \Pol_V(\Lambda)$ as a subring of $\cP_c(V, \bQ)$ and we define $\cP_c(V, \bZ)$ to be the smallest subring containing $C_c(\Lambda_{\bhatZ}, \bZ)\otimes \Pol_V(\Lambda)$ for all $\Lambda$, i.e.\
\begin{equation*}
\cP_c(V, \bZ)\,\, =\,\, \Image\left(\bigoplus_{\Lambda} C_c(\Lambda_{\bhatZ}, \bZ)\otimes \Pol_V(\Lambda)\to \cP_c(V, \bQ)\right).
\end{equation*}
It is a $\GL(V)$-stable subring of $\cP_c(V, \bQ)$. We remark that for $h\in \cP_c(V, \bZ)$ we have $h(v)\in \bZ$ for all $v\in V$. Note also that $\cP_c(V, \bQ) = \cP_c(V, \bZ)_{\bQ}$. For an arbitrary ring $R$ we put $\cP_c(V, R) = \cP_c(V, \bZ)_R$. 

Now we turn to  the Solomon-Hu pairing. Let $\sL = \{\sL_i\}_{i\in I} = \{(\ell_i, L_i, L_i')\}_{i\in I}$ be a collection of triples as above involving a set of prime elements $T = \{\ell_i\mid i\in I\}$ such that $\ell_{i} \ne \ell_{i'}$ for $i\ne i'$ ($T$ consisting of a single prime $\ell$ is allowed). We define \[ \cP_c(V, \bQ)^T = C_c(V_{\bhatZ^T}, \bZ)\otimes \Pol_V\] and more generally for a $\bQ$-algebra $R$ we put $\cP_c(V, R)^T = C_c(V_{\bhatZ^T}, R)\otimes \Pol_V$. For $f\in \cN^{\sL}$, $\phi\in C_c(V_{\bhatZ^T}, \bZ)$ and $P\in \Pol_V$ we can evaluate $P$ at $\llangle f, \phi \rrangle^{\sL}\in \bZ\{\!\{V\}\!\}^{\sL}\subseteq \bQ[\![V]\!]$ (i.e.\ we can apply \eqref{pol2} to $\llangle f, \phi \rrangle^{\sL}$). In this way we obtain for any $\bQ$-algebra $R$ an $\Gamma^{\sL}$-equivariant pairing 
\begin{equation}
\label{solhupol1}
\llangle\,\,\,, \,\,\, \rrangle\,:\,\, \cN^{\sL} \times\cP_c(V, R)^T  \, \lra\, R,\,\, \llangle f, \phi\otimes P \rrangle \, = \, P(\llangle f, \phi \rrangle^{\sL}).
\end{equation}
Let $\cP_c(V, \bZ)^\sL$ denote the smallest subring of $\cP_c(V, \bQ)^T$ containing all $C_c(\Lambda_{\bhatZ^T}, \bZ)\otimes \Pol_V(\Lambda)$ where $\Lambda$ runs through all lattices of $V$ with $\Lambda_{\ell_i} = L_i$ for all $i\in I$. The $\Gamma^{\sL}$-equivariant embedding
\begin{equation*}
\cP_c(V, \bQ)^T\,\lra\,\cP_c(V, \bQ),  \,\, \phi\otimes P\mapsto (\phi\otimes \bigotimes_{i\in I} \phi_{\sL_i})\otimes P
\end{equation*}
maps $\cP_c(V, \bZ)^\sL$ into $\cP_c(V, \bZ)$. 

\begin{lemma} 
\label{lemma:eisenint}
For $f\in \cN^{\sL}$ and $h\in \cP_c(V, \bZ)^\sL$ we have $\llangle f, h \rrangle\in \bZ$ (resp.\ $\llangle f, h \rrangle\in \bZ[1/\ell])$ if $T$ has at least two elements (resp.\ if $T = \{\ell\}$).
\end{lemma}

{\em Proof.} We assume that $T$ contains at least two primes and leave the case $T=\{\ell\}$ to the reader. Let $\Lambda\subseteq V$ be lattice with $\Lambda_{\ell_i} = L_i$ for all $i\in I$ and let $\phi \in C_c(\Lambda_{\bhatZ^T}, \bZ)$ and $P\in \cP(\Lambda)$. Since $\llangle f, \phi \rrangle^{\sL}\in \bZ[\![\Lambda]\!]$ we have $\llangle f, \phi\otimes P \rrangle = P(\llangle f, \phi \rrangle^{\sL})\in \bZ$. In fact we can approximate $\llangle f, \phi \rrangle^{\sL}$ by an element of the group ring $\bZ[\Lambda]$ modulo any power of the augmentation ideal and by assumption $P$ maps $\bZ[\Lambda]$ to $\bZ$. \enddemo

Thus if $\sharp(T)\ge 2$ by restricting \eqref{solhupol1} to the submodule $\cP_c(V, \bZ)^\sL$ of $\cP_c(V, \bQ)^T$ we obtain a $\Gamma^{\sL}$-equivariant pairing 
\begin{equation}
\label{solhupol2}
\llangle\,\,\,, \,\,\, \rrangle\,:\,\, \cN^{\sL} \times\cP_c(V, \bZ)^{\sL}  \, \lra\, \bZ
\end{equation}
More generally if $R$ is a ring we put $\cP_c(V, R)^\sL = \cP_c(V, \bZ)^\sL\otimes R$ so we obtain a pairing 
\begin{equation}
\label{solhupol2a}
\llangle\,\,\,, \,\,\, \rrangle\,:\,\, \cN^{\sL} \times\cP_c(V, R)^{\sL}  \, \lra\, R
\end{equation}
Note that $\cP_c(V, R)^{\sL} = \cP_c(V, R)^T$ if $R$ is a $\bQ$-algebra.

Analogously to \eqref{solhumod} the pairings \eqref{solhu4b}, \eqref{solhupol2} induces homomorphism of $\Gamma^{\sL}$-modules 
\begin{eqnarray}
\label{solhumodint1}
&& \Hom(\sR,\cN^{\sL})(\sgn)\, \lra \, \sM(C_c(V_{\bhatZ^T}, \bZ), \bZ)(\sgn),\\
\label{solhumodint2}
&& \Hom(\sR,\cN^{\sL})(\sgn) \, \lra \,  \sM(\cP_c(V, \bZ)^{\sL}, \bZ)(\sgn).
\end{eqnarray}

\begin{definition}
\label{definition:eisen}
Let $\sL = \{\sL_i\}_{i\in I}=\{(\ell_i, L_i, L_i')\}_{i\in I}$ be a collection of triples as before involving a finite set of prime numbers $T= \{\ell_i\mid i\in I\}$.

\noi (a) Assume that $T$ contains at least two elements. The integral Eisenstein cocycle
\[ 
\Eis_{\sL}\,\in H^{n-1}(\Gamma^{\sL},\sM(\cP_c(V, \bZ)^{\sL}, \bZ)(\sgn))
\]
is defined as the image of $\Xi_{\Sh,\sL}$ under the map induced by \eqref{solhumodint2}. 

\noi (b) Assume either $\sharp{T}\ge 2$ or $T=\{\ell\}$ and $\ell\ge n+2$. The truncated integral Eisenstein cocycle 
\[
\Eis_{\sL}^0\,\in H^{n-1}(\Gamma^{\sL},\sM(C_c(V_{\bhatZ^T}, \bZ), \bZ)(\sgn))
\] 
is defined as the image of $\Xi_{\Sh,\sL}$ under the map induced by \eqref{solhumodint1}. 
\end{definition}

More generally by using the pairing \eqref{solhupol2a} for an arbitrary ring $R$ we can define the $R$-valued Eisenstein cocycle 
\[
\Eis_{\sL,R}\,\in H^{n-1}(\Gamma^{\sL},\sM(\cP_c(V, R)^{\sL}, R)(\sgn)).
\]

\section{Specialization to totally real fields}
\label{section:fields}

\subsection{Integral values of Shintani zeta function}

We keep the notation of the previous section and fix a basis $j_1, \ldots, j_n: V_{\bR}\to \bR$ of the dual of $V_{\bR}$. We assume that the restriction of $j_\nu$ to $V$ has trivial kernel for all $\nu=1,\ldots, n$. A vector $v\in V_{\bR}$ is called totally positive (with respect to $j_1, \ldots, j_n$) if $j_{\nu}(v) > 0$ for all $\nu=1,\ldots, n$. By $V_{\bR,+}$ we denote the set of totally positive vectors. A rational open cone $C=C(v_1,  \ldots, v_m)$ contained in $V_{\bR,+}$ will be called positive. We denote by $\cK_+$ the abelian group of functions $V_{\bR} \to \bZ$ generated by characteristic functions of positive rational open cones. Note that $\cK^+\cap \cL=\{0\}$ so that we can view $\cK_+$ as a subgroup of $\cN=\cK/\cL$. 

Given $f\in \cK_+$ and $\Phi \in C_c(V_{\bhatZ}, \bZ)$ we consider the Dirichlet series 
\begin{equation}
\label{shintanizeta}
L(f, \Phi; s) \,=\, \sum_{v\in V} \, f(v)\,\Phi(v)\,\Norm(v)^{-s}
\end{equation}
where $\Norm(v) = \prod_{\nu=1}^n j_\nu(v)$. It is known to converge absolutely for $\Real(s) >1$ and extend to the whole complex plane holomorphically except for a simple pole at $s=0$. Moreover if $C$ and $\Phi$ are as in Prop. \ref{prop:eisenint} then $L(1_C,\Phi,s)$ is holomorphic and its values at $s= 0, -1, -2, \ldots$ can be expressed in terms of the pairing \eqref{solhupol1}. In fact we have

\begin{lemma} 
\label{lemma:shintani}
Let $\sL= (\ell, L, L')$ be a triple as in section \ref{subsection:cassou}, let $C=C(v_1,  \ldots, v_m)$ with $v_i\in (V\cap L)-L'$ and let $\Phi=\phi\otimes \phi_{\sL}$ with $\phi\in C_c(V_{\bhatZ^{\ell}}, \bZ)$ and $\phi_{\sL_i}\in C_c(V_{\ell},\bZ)$ given by \eqref{cassousmooth}. Then the function $L(1_C, \Phi; s)$ extends holomorphically to the whole complex plane. Its values at integers $s=k$, $k\in \bZ_{\le 0}$ are given by
\begin{equation*}
L(1_C, \Phi ;k)\,\, =\,\, \llangle 1_C, \Phi\otimes \Norm^{-k}\rrangle.
\end{equation*}
\end{lemma}

{\em Proof.} With the notation and assumption as in the proof of Prop.\ \ref{prop:eisenint} and an appropriate choice of a lattice $\Lambda$ with $\Lambda_{\ell}=L$ and $\supp(\Phi_0)\subseteq \Lambda_{\bhatZ}$ it follows from \eqref{eisenint2}
\begin{equation*}
\llangle 1_C, \Phi\otimes \Norm^k\rrangle\,=  \, -\sum_{\psi} \sum_{v\in \Lambda\cap P}  \psi(v)\Phi_0(v) \frac{\Norm(v)^k}{\prod_{i=1}^m(1- \psi(v_i)\Norm(v_i)^k)}.
\end{equation*} 
On the other hand for $L(1_C, \Phi; s)$ we have
\begin{eqnarray*}
&& L(1_C, \Phi ; s) \,\,=\,\,  - \sum_{\psi}\sum_{v\in \Lambda\cap C}  \psi(v) \Phi_0(v) \Norm(v)^{-s}\\ 
&& =\, - \sum_{\psi} \sum_{v\in \Lambda \cap P} \psi(v) \Phi_0(v) 
 \sum_{n_1,\ldots, n_m=0}^{\infty} \psi\!\left(\sum_{i=1}^m n_i v_i\right)\Norm(v + \sum_{i=1}^m n_i v_i)^{-s}.
\end{eqnarray*}
Thus it is enough to show that for every unitary character $\psi: \Lambda\to \bC^*$ with $\psi(v_i) \ne 1$ for all $i=1, \ldots, m$ the Dirichlet series 
\[
 \sum_{n_1,\ldots, n_m=0}^{\infty} \psi\!\left(\sum_{i=1}^m n_i v_i\right) \Norm(v + \sum_{i=1}^m n_i v_i)^{-s}
\]
extends to a holomorphic function on the whole complex plane and that  value at $s=k$ is equal to $\frac{\Norm(v)^{-k}}{\prod_{i=1}^m\, (1 -  \psi(v_i)\Norm(v_i)^{-k})}$. This is well-known (see \cite{shintani}, Prop.\ 9). 
\enddemo

\subsection{Zeta functions of totally real fields}
\label{subsection:totalreal}

We fix now a totally real number field $F$ of degree $n$ over $\bQ$ and let $K/F$ be a finite abelian extension with Galois group $G$. Let $S$ be a finite set of nonarchimedean places of $F$ containing all places which are ramified in $K$. For $\sigma\in G$ we recall the definition of the partial zeta function  
\begin{equation*}
\zeta_S(\sigma,s)\,\, =\,\,\sum_{(\fa,S)=1, \sigma_{\fa}=\sigma}\, \Norm(\fa)^{-s}.
\end{equation*} 
Here the sum is taken over all ideals $\fa\subseteq \cO_F$ that are relatively prime to the elements in $S$ and such that their image $\sigma_{\fa}\in G$ under the Artin map is equal to $\sigma$. The series is absolutely convergent for $\Real(s) > 1$ and extends to a meromorphic function on the whole complex plane with a single pole at $s=1$.

By a theorem of Siegel and Klingen the values at nonpositive integers $s=k$ of $\zeta_S(\sigma,s)$ are rational numbers. They can be expressed in terms of the Eisenstein cocycle \cite{sczech, cdg}. In this section we  describe this relation within our framework (as it differs from \cite{cdg}). We choose an auxiliary finite non-empty set $T$ of primes of $F$ disjoint from $S$ and define the $\bC[G]$-valued zeta function $\zeta_{S, T}(K/F, s)$ by the following identity
\begin{equation*}
\zeta_{S, T}(K/F, s) \,\,=\,\,  \prod_{\fq\in T} (1-\Norm(\fq)^{1-s}[\sigma_{\fq}^{-1}]) \, \sum_{\sigma\in G}\zeta_S(\sigma,s)[\sigma^{-1}].
\end{equation*} 
The partial zeta functions $\zeta_{S,T}(\sigma,s)$ are defined as the different components of $\zeta_{S, T}(K/F, s)$, i.e.\ we have the following identity in $\bC[G]$
\begin{equation*}
\zeta_{S, T}(K/F, s) \,\, =\,\, \sum_{\sigma\in G}\zeta_{S,T}(\sigma,s)[\sigma^{-1}].
\end{equation*}
A character $\chi\colon G\to \bC^*$ induces a ringhomomorphism $\chi\colon\bC[G]\to \bC$ and we obtain the usual $L$-functions by $L_{S, T}(\chi, s) =\chi^{-1}(\zeta_{S, T}(K/F, s))$.

We denote by $T_{\bQ}= \{q \mid q\in \fq \, \mbox{for some $\fq\in T$}\,\}$ the set of prime numbers that underlie $T$ and by $\barT$ the set of places of $F$ that lie above a prime in $T_{\bQ}$. We consider the following conditions on $T$.
\begin{itemize}
\item[\bf(A1) ] Over every $q\in T_{\bQ}$ lies exactly one $\fq\in T$ and we have $\Norm(\fq) = q$.
\item[\bf(A2) ] $\barT$ and $S$ are disjoint.
\item[\bf(A3) ] $T$ has at least two elements.
\item[\bf(A3a)] $T$ has at least two elements or $T= \{\fq\}$ and $q=\Norm(\fq)\ge n+2$.
\end{itemize}

In order to express the relation between the value of $\zeta_{S, T}(K/F, s)$ at $k\in \bZ_{\le 0}$ we introduce some notation. We choose an isomorphism of $\bQ$-vector spaces $F\cong V$ which allows us to identify $F$ with $V$ in the following (so the reader may think of $V=F$ considered just as a $\bQ$-vector space). Therefore we can view the torus $\Res_{F/\bQ}\bG_m$ as an algebraic subgroup of $\GL(V)$. Also we associate to every prime $\fq$ with underlying prime number $q$ the triple $\sL_{\fq} = (\, \cO_F\otimes \bZ_q, \fq\otimes \bZ_q)$. Thus we have a collection of triples $\sL = \{\sL_{\fq}\}_{\fq\in T}$ as in section \ref{section:eisenstein}. We also identify $C_c(\bA_F^{\infty}, \bZ)$, $C_c(\bA_F^{\barT,\infty}, \bZ)$ etc.\ with $C_c(V_{\bhatZ}, \bZ)$, $C_c(V_{\bhatZ^T}, \bZ)$ etc. (here and in the following we write $\bhatZ^T$ for $\bhatZ^{T_{\bQ}}$). Note that $F^{\barT}\subseteq \Gamma^{\sL} \cap F^*$.

The isomorphism $F\cong V$ provides $V$ with the structure of an $F$-vector space hence induces an $F^*$-action on $\sR$ (recall that $\sR$ was defined in \S\ref{s:shintanicocycle}). For $v\in S_{\infty}$ we denote by $\sigma_v: F\to F_v=\bR$ the corresponding embedding and by $\sgn_v:F_v^* \to\{\pm 1\}$ the sign character. There exists an 
unique $\{\pm 1\}$-orbit $\sF_v\subseteq \sR$ such that 
\begin{equation*}
x \cdot Q \,\, =\,\, \sgn_v(x) \,Q
\end{equation*}
for all $x\in F^*$ and $Q\in \sF_v$. In fact if $\{e_v\}_{v\in S_{\infty}}$ is a basis of $V_{\bR}$ such that $x \cdot e_v = \sigma_v(x) e_v$ for all $v\in S_{\infty}$ and $x\in F^*$ then $\sF_v$ consists of the two rays $\pm \bR_+ e_v$.
The embedding $\sigma_v: F\cong V\to \bR$ can be extended to a linear form $j_v: V_{\bR} \cong F_{\bR} \to \bR$ and we use the set $\{j_v\}_{v\in S_{\infty}}$ to define the notion of totally positive elements in $V_{\bR}$. Note that the restriction of $\Norm = \prod_{v} j_v$ to $F\cong V$ is the usual norm $\Norm_{F/\bQ}$ of the extension $F/\bQ$.

For $v\in S_{\infty}$, $Q\in \sF_v$, a $\Gamma^{\sL}$-module $A$ and abelian group $B$,  the ``evaluation at $Q$" map 
\[ \ev_{Q}\colon \sM(A, B) \to  \Hom(A, B),\qquad \beta \mapsto  \beta(Q, \wcdot) \] is $F^{\barT, v}$-equivariant. Hence the pair $(F^{\barT, v}\hookrightarrow \Gamma^{\sL}, \ev_Q)$ induces a homomorphism
\begin{equation}
\label{eval1}
H^{n-1}(\Gamma^{\sL}, \sM(A, B)(\sgn)) \to H^{n-1}(F^{\barT, v}, \Hom(A, B)(\sgn)).
\end{equation}

\begin{definition}
\label{definition:eisenF}
(a) Assume that (A1), (A2) and (A3) hold. We define the cohomology class
\begin{equation*} 
\Eis_{F, T, Q} \,\, \in H^{n-1}(F^{\barT,v}, \Hom(\cP_c(V, \bZ)^{\sL}, \bZ)(\sgn))
\end{equation*}
as the image of the Eisenstein cocycle $\Eis_{\sL}$ under the homomorphism \eqref{eval1} for $A= \cP_c(V, \bZ)^{\sL}$.

\noi (b) Assume that (A1), (A2) and (A3a) hold. The cohomology class
\begin{equation*} 
\Eis_{F, T, Q}^0 \,\, \in H^{n-1}(F^{\barT,v}, \Hom(C_c(V_{\bhatZ^T}, \bZ), \bZ)(\sgn)) 
\end{equation*} 
is defined as the image of $\Eis_{\sL}^0$ under the homomorphism \eqref{eval1} for $A= C_c(V_{\bhatZ^T}, \bZ)$.
\end{definition}

In order to define the Eisenstein cocycle on the full group $F^{\barT}$ rather than the subgroup $F^{\barT, v}$, we note that $\ev_Q + \ev_{-Q}$ is $F^{\barT}$-equivariant and hence induces a homomorphism
\begin{equation}\label{evalplusminus}
H^{n-1}(\Gamma^{\sL}, \sM(A, B)(\sgn)) \to H^{n-1}(F^{\barT}, \Hom(A, B)(\sgn)).
\end{equation}
Parallel to Definition \ref{definition:eisenF} above, we define
\begin{align} 
\Eis_{F, T, v} \,\, &\in H^{n-1}(F^{\barT}, \Hom(\cP_c(V, \bZ)^{\sL}, \bZ)(\sgn))  \label{e:defeisv} \\
\Eis_{F, T, v}^0 \,\, & \in H^{n-1}(F^{\barT}, \Hom(C_c(V_{\bhatZ^T}, \bZ), \bZ)(\sgn)) \nonumber
\end{align}
as the image of $\Eis_{\sL}$ under (\ref{evalplusminus}) in the cases (a) and (b), respectively.

If $H$ is a subgroup of $F^{\barT,v}$ then (by abuse of notation) we shall denote the restriction of $\Eis_{F, T, Q}$ and $\Eis_{F, T, Q}^0$ to $H$-cohomology also simply by $\Eis_{F, T, Q}$. In fact for most of our applications it is enough to consider the restriction of $\Eis_{F, T, Q}^0$ to the group of totally positive elements in $F^{\barT}$, i.e.\ we can work with the class 
\begin{equation}
\label{eisennosigntruc}
\Eis_{F, T, Q}^0 \,\, \in H^{n-1}(F^{\barT}_+, \Hom(C_c(V_{\bhatZ^T}, \bZ), \bZ)).
\end{equation}

\subsection{The homomorphism $\partial$}
Let $A$ be a locally profinite abelian group. We construct an $(\bA_F^{\barT,\infty})^*$-equivariant map
\begin{equation}
\label{delta3}
\Delta_S^{\barT}: \cCcd(S, A)^{\barT} \,\lra \, \Ccd(\bA_F^{\barT,\infty}, A) \cong \Ccd(V_{\bhatZ^T}, A).
\end{equation}
as follows. There exist canonical homomorphisms
\begin{eqnarray*}
&& \Ccd(\prod_{v\in S} F_v^*, A) \otimes  \cC_c(\emptyset, \bZ)^{S\cup \barT,\infty}\, \lra \, \cCcd(S, A)^{\barT,\infty}\\
&& \Ccd(\prod_{v\in S} F_v, A) \otimes C_c(\bA_F^{S\cup \barT,\infty}, \bZ)\,\lra \, \Ccd(\bA_F^{\barT,\infty}, A)
\end{eqnarray*}
(see \eqref{multfunctions4}, \eqref{multfunctions5}). In fact the first map is an isomorphism.
Since $(\bA^{S\cup \barT, \infty}_F)^*/U^{S\cup \barT, \infty}$ is isomorphic to the group of fractional ideals $\cI^{S\cup \barT}$ of $F$ that are relatively prime to $S\cup \barT$, the ring $\cC_c^0(\emptyset, \bZ)^{S\cup \barT, \infty}$ can be identified with the group ring $\bZ[\cI^{S\cup \barT}]$.
We define \eqref{delta3} as the tensor product $\Delta_S^{\barT} = \delta_S \otimes  \delta^{S\cup \barT}$ where $\delta_S\colon \Ccd(\prod_{v\in S} F_v^*, A) \to \Ccd(\prod_{v\in S} F_v, A)$ is the inclusion \eqref{extensionsbyzero} and $\delta^{S\cup \barT}$ maps a fractional ideal $\fa\in \cI^{S\cup \barT}$ to the characteristic function of $\,\widehat{\fa}^{S\cup \barT} =  \fa\, (\prod_{\fp\not\in S\cup \barT}\cO_{\fp})$.

Given a decomposition $S_{\infty}= S_3\cup S_4$ into disjoint subsets consider the composition
\begin{align}
\label{delta6a}
& \Delta_*: H_{n-1}(F^*, \cCcd(S, A)) \stackrel{\eqref{augep2}_*}{\lra}\, H_{n-1}(F^*, \cCcd(S, A)^{S_3}(\sgn_{S_3})) \nonumber \\
& \hspace{1cm}\stackrel{(*)}{\cong} \, H_{n-1}(F^{\barT\cup S_4}, \cCcd(S, A)^{\barT, \infty}(\sgn)) \\
& \hspace{1cm}\stackrel{\eqref{delta3}_*}{\lra} \,  H_{n-1}(F^{\barT\cup S_4}, \Ccd(V_{\bhatZ^T}, A)(\sgn)).\nonumber
\end{align}
For the isomorphism $(*)$ note that $\cCcd(S, A)^{S_3} = \Ind_{F^{\barT\cup S_4}}^{F^*}  \cCcd(S, A)^{\barT, \infty}$. We define 
\begin{equation}
\label{delta6}
 \partial\,=\,\partial_{S, S_3, S_4}\colon H^0(F^*, \cCd(S, A))\to H_{n-1}(F^{\barT\cup S_4}, \Ccd(V_{\bhatZ^T}, A)(\sgn)).
\end{equation}
as the composition of \eqref{delta2} with \eqref{delta6a}.

For later use we remark that for a decomposition of $S$ into disjoint subsets $S= S_1\cup S_2$ the map \eqref{delta6a} can be defined for the larger group of coefficients $\cCcd(S_1, S_2, A)$ as well, i.e.\ we obtain a map 
\begin{equation}
\label{delta6b}
\Delta_*: H_{n-1}(F^*, \cCcd(S_1, S_2, A)) \, \lra \,  H_{n-1}(F^{\barT\cup S_4}, \Ccd(V_{\bhatZ^T}, A)(\sgn))
\end{equation}
by replacing $\cCcd(S, A)$ by $\cCcd(S_1, S_2, A)$ everywhere in the definition of \eqref{delta6a} above. 
For that we only need to remark that the map \eqref{delta3} can be extended to $\cCcd(S_1, S_2, A)^{\barT}$ by replacing $\delta_S$ with the inculsion $\Ccd(\prod_{v\in S_1} F_v\times \prod_{v\in S_2} F_v^*, A) \stackrel{\eqref{extensionsbyzero}}{\lra} \Ccd(\prod_{v\in S} F_v, A).$

We consider now the case where $A=R$ is a ring equipped with the discrete topology so that \eqref{delta6}  is the map
\begin{equation*}
\partial\colon H^0(F^*, \cC^0(S, R))\,\lra\, H_{n-1}(F^{\barT\cup S_4}, C_c^0(V_{\bhatZ^T}, R)(\sgn)).
\end{equation*}
We need certain ``$k$-twisted" variants of $\partial$, i.e.\ we define maps
\begin{equation}
\label{delta7}
\partial^k=\partial^k_{S, S_3, S_4}: H^0(F^*, \cC^0(S, R))\to H_{n-1}(F^{\barT\cup S_4}, \cP_c(V, R)^{\sL}(\sgn)).
\end{equation}
for every $k\in \bZ_{\le 0}$.
For this we first introduce ``$k$-twisted" versions $\Delta(k)_S^{\barT}$ of \eqref{delta3} for $k\in \bZ_{\le 0}$. Define the ``idele norm character" $\bN\colon \bA_F^*\to \bQ^*, (x_v)_v\mapsto \prod_v \bN_v(x_v)$ by
\begin{equation*}
\bN_v(x) \,=\,  \left\{\begin{array}{ll} |x|_v^{-1} &  \mbox{if $v$ is nonarchimedean,}\\
\sgn(x) &  \mbox{if $v$ is archimedean.}
\end{array}\right.
\end{equation*}
Note that for $a\in F^*$ we have $\bN(a)=\Norm_{F/\bQ}(a)=\colon \Norm(a) $. Hence
\begin{equation}
\label{delta4}
\Delta(k)_S^{\barT}: \cC_c^0(S, \bZ)^{\barT, \infty}(\sgn^k)  \,\lra \, \cP_c(V, \bQ)^T, \qquad \phi\mapsto \Delta_S^{\barT}(\bN^k\cdot \phi) \otimes \Norm^{-k}
\end{equation}
is $F^{\barT}$-equivariant (recall that $\sgn$ denotes the character \eqref{sign}; we view it as a character of $F^*$ via the embedding $F^*\hookrightarrow \GL(V)$). Here $\bN^k\cdot \phi$ denotes the function 
\[ (\bA_F^{\barT, \infty})^*\to \bQ, \qquad x=(x_v)_{v\not\in \barT, v\nmid \infty} \mapsto \prod_{v\not\in \barT, v\nmid \infty} \bN_v^k(x_v) \phi(x). \] 

\begin{lemma} 
The image of the map \eqref{delta4} is contained in $\cP_c(V,\bZ)^{\sL}$. 
\end{lemma}

{\em Proof.} As an $F^{\barT}$-module $\cC_c(S, \bZ)^{\barT,\infty}$ is generated by functions of the type
$1_{U'}\otimes 1_{aU^{S\cup \barT, \infty}}$ where $U'$ is an open subgroup of $U_S$ and $a = (a_v)\in (\bA_F^{S\cup \barT,\infty})^*$. Therefore it suffices to show that $\Delta_S^{\barT}(1_{U'}\otimes 1_{aU^{S\cup \barT,\infty}})$ is contained in $\cP_c(V,\bZ)^{\sL}$. Let $\fa\in \cI^{S\cup \barT}$ be the fractional ideal corresponding to $a$. We have 
\begin{equation}
\label{generator2}
\Delta(k)_S^{\barT}(1_{U'}\otimes 1_{aU^{S\cup \barT,\infty}}) \, =\, (1_{U'}\otimes 1_{\widehat{\fa}^{T\cup S}})\otimes (\Norm(\fa)^k\Norm^{-k})
\end{equation}
If $\Lambda\subseteq V$ corresponds to $\fa$ under the identification $V\cong F$ then $(\Norm(\fa)^k\Norm^{-k})\in \Pol_V(\Lambda)$ and therefore \eqref{generator2} is contained in $C_c(\Lambda_{\bhatZ^T} ,\bZ)\otimes \Pol_V(\Lambda)\subseteq \cP_c(V,\bZ)^{\sL}$.\enddemo

The Lemma shows that \eqref{delta4} can be viewed as a map $\cC_c(S, \bZ)^{\barT,\infty}\to \cP_c(V,\bZ)^{\sL}$. Hence for any ring $R$ we can define 
\begin{equation}
\label{delta4a}
\Delta(k)_S^{\barT}: \cC_c^0(S, R)^{\barT,\infty}(\sgn^{k+1})  \,\lra \, \cP_c(V, R)^{\sL}(\sgn)
\end{equation}
by taking the tensor product of \eqref{delta4} with $\id_R$. We define \eqref{delta7} as the composite map
\begin{eqnarray*}
&& \partial^k\,=\,\partial^k_{S, S_3, S_4}\colon H^0(F^*, \cC^0(S, R)) \,\stackrel{\eqref{delta2}}{\lra}\, H_{n-1}(F^*, \cC_c^0(S, R))\\
&&\stackrel{\eqref{augep2}_*}{\lra}\, H_{n-1}(F^*, \cC_c^0(S, R)^{S_3}(\sgn_{S_3}^{k+1}))\,\stackrel{(*)}{\cong} \, H_{n-1}(F^{\barT\cup S_4}, \cC_c^0(S, R)^{\barT, \infty}(\sgn^{k+1}))\nonumber\\
&& \hspace{2cm}\stackrel{\eqref{delta4a}_*}{\lra} \,  H_{n-1}(F^{\barT\cup S_4}, \cP_c(V, R)^{\sL}(\sgn)).\nonumber
\end{eqnarray*}

For $v\in S_{\infty}$ consider the cap-product pairing 
\begin{equation*}
H^{n-1}(F^{\barT, v}, \Hom(\cP_c(V, \bZ)^{\sL}, \bZ)(\sgn)) \times H_{n-1}(F^{\barT, v}, \cP_c(V, R)^{\sL}(\sgn))\to R
 \end{equation*}
induced by the canonical map $\Hom(\cP_c(V, \bZ)^{\sL}, \bZ) \times \cP_c(V, R)^{\sL} \to R$. Using \eqref{delta6} for $S_3=S_{\infty}-\{v\}$ and $S_4= \{v\}$ we obtain a pairing
\begin{equation*}
H^{n-1}(F^{\barT, v}, \Hom(\cP_c(V, \bZ)^{\sL}, \bZ)(\sgn)) \times H^0(F^*, \cC^0(S, R))\to R
 \end{equation*}
given by $(x,y)\mapsto x\cap \partial^k(y)$.   Next we apply this pairing with the Eisenstein cocycle to deduce
integrality results on partial zeta functions of totally real fields.

\subsection{Stickelberger elements and Shintani zeta functions}
Let
\begin{equation*}
\rec_{K/F}:\bA_F^*\,\lra\, G\subseteq \bZ[G]^*,\, x= (x_v)_v\mapsto\rec(x)= \prod_v (x_v, K/F)_v
\end{equation*}
be the global reciprocity map. Here $(x, K/F)_v$ denotes the local norm residue symbol for every place $v$ of $F$ and $x\in F_v^*$. We can view $\rec_{K/F}$ as an element of $H^0(F^*, \cC^0(S, \bZ[G]))$ so we can make the following definition.

\begin{definition}
\label{definition:stickelberger}
Assume that (A1), (A2) and (A3) hold. For $k\in \bZ_{\le 0}$ and $Q\in \sF_v$ the Stickelberger element is defined by 
\begin{equation*}
\Theta_{S, T}(K/F, k)\,\, = \,\, \Eis_{F, T, Q}\cap\partial^k(\rec_{K/F})\in \bZ[G].
\end{equation*}
If only (A3a) holds instead of (A3) then we define $\Theta_{S, T}(K/F, 0)$ by
\begin{equation*}
\Theta_{S, T}(K/F, 0)\,\, = \,\, \Eis_{F, T, Q}^0\cap\partial(\rec_{K/F})\in \bZ[G].
\end{equation*}
\end{definition}

\begin{remark}
\label{remark:stickelberger}
\rm \noi (a) Using the fact that $\partial^k_{S, S_3, S_4}= \cor\circ \partial^k_{S, S_3', S_4'}$ if $S_4 \subseteq S_4'$ (where $\cor$ denotes the corestriction from $F^{\barT\cup S_4'}$- to $F^{\barT\cup S_4}$-homology) one shows that $\Theta_{S, T}(K/F, k)$ does not depend on the choice of the decomposition of $S_{\infty} = S_3 \cup S_4$. The choice of $S_3$, $S_4$ will be important in the proof of Thm.\ \ref{theo:highervanishing} however.

\noi (b) Let $L/F$ be an abelian extension with $L\supseteq K$ and put $\widetilde{G} =\Gal(L/F)$.  Assume that $L/F$ is unramified outside $S$ and let $\pi\colon \bZ[\widetilde{G}]\to \bZ[G]$ be the canonical projection. Since $\pi\circ \rec_{L/F} = \rec_{K/F}$ we have  
\[
\pi(\Theta_{S, T}(L/F, k)) = \Theta_{S, T}(K/F, k).
\]
\noi (c) If $v$ is a nonarchimedean place of $F$ with $v\not\in S\cup \barT$ and  we put $S'= S\cup \{v\}$ then we have
\begin{equation}
\label{changeS0}
\Theta_{S', T}(K/F, k) = (1-\Norm(v)^{-k}[\sigma_v^{-1}])\Theta_{S, T}(K/F, k). 
\end{equation}
In fact we have 
\begin{equation*}
\partial^k_{S'}(\rec_{K/F}) \, =\, (1-\Norm(v)^{-k}[\sigma_v^{-1}])\partial^k_{S}(\rec_{K/F})
\end{equation*}
which can be easily deduced from the commutativity of the diagram
\begin{equation}
\label{changeS2}
\begin{CD}
\cC_c(S, \bZ)^{\barT,\infty} @> \phi \mapsto \phi- \Norm(v)^{-k}[\varpi_v] \phi>> \cC_c(S, \bZ)^{\barT,\infty}\\
@VV \incl V @VV \Delta(k)_S^{\barT} V\\
\cC_c(S', \bZ)^{\barT,\infty} @> \Delta(k)_{S'}^{\barT} >> \cP_c(V, \bZ)^{\sL}
\end{CD}
\end{equation}
and the fact that $[\varpi_v] \rec_{K/F}= [\sigma_v^{-1}]\rec_{K/F}$.
Here $[\varpi_v]\in (\bA_F^{\barT})^*$ denotes the idele whose component at $v$ is the uniformizer $\varpi_v\in F_v$ and whose other components are $=1$.
\end{remark}

Of course (b) and (c) also follow  immediately from:

\begin{prop} 
\label{prop:zetaelement}
For $k\in \bZ_{\le 0}$ we have $\Theta_{S, T}(K/F, k) = \zeta_{S, T}(K/F, k)$. 
\end{prop}

\begin{remark} 
\rm In particular we obtain $\zeta_{S, T}(K/F, k)\in \bZ[G]$ and also that $\Theta_{S, T}(K/F, k)$ does not depend on the choices of $v$ and $Q$.
\end{remark}

In order to prove Prop.\ \ref{prop:zetaelement} we first recall how to express $\zeta_{S, T}(K/F, s)$ in terms of Shintani zeta functions. The definition \eqref{shintanizeta} can be extended naturally to functions $\Phi\in \cP_c(V_{\bhatZ}, \bZ[G])$ (in this case $L(f, \Phi, s)$ has values in $\bC[G]$). For a fractional ideal $\fa\in \cI^{S\cup T}$ we consider the function $\phi_{\fa} \in \cP_c(V_{\bhatZ^T}, \bZ[G])$ defined by
\begin{equation*}
\phi_{\fa}(x) \, = \, [\rec_S(x)\sigma_{\fa}] (1_{\widehat{\fa}^{T\cup S}}\otimes 1_{U_S})(x)
\end{equation*}
where $\rec_S\colon \bbI^T\to G$ is given by $\rec_S((x_v)_{v\not\in T}) = \prod_{v\in S}(x_v,K/F)_v\in G$ i.e.\ $\rec_S((x_v)_{v\not\in T})$ is the product of the local norm residue symbols $(x_v,K/F)_v$ at places in $S$. Note that for $y\in E_{S\cup T,+}$ we have $\phi_{y\fa}\, =\, y\phi_{\fa}$ because of $\rec_S(y)\sigma_{(y)} = 1$. We also define $\Phi_{\fa} \in \cP_c(V_{\bhatZ}, \bZ[G])$ by 
\begin{equation*}
\Phi_{\fa} \, = \, \phi_{\fa}\otimes \bigotimes_{\fq\in T} \phi_{\sL_{\fq}}.
\end{equation*}
Let $\cA\subseteq V_{\bR,+}$ be a Shintani decomposion for $F$. Recall that this means that $\cA$ can written as a finite disjoint union of rational open cones and that $V_{\bR,+} \,\,\, =\,\,\, \overset{\cdot}{\bigcup}_{\epsilon \in E_+} \ep \cA$. We fix represenatives $\fa_1, \ldots, \fa_h$ of the different ideal classes in the narrow class group $\Cl^+(F)$ that are relatively prime to $S\cup \barT$. 

\begin{lemma} 
\label{lemma:shintani2}
We have $\zeta_{S, T}(K/F, s)\, =\,\sum_{i=1}^h\, \Norm(\fa_i)^s L(1_{\cA}, \Phi_{\fa_i}, s)$.
\end{lemma}

{\em Proof.} This is a well-known computation (see e.g.\ \cite{spiess2}). Let $\fa\in \cI^S$ and let $\fA$ be its class in $\Cl^+(F)$. The map $x\mapsto x\fa^{-1}$ induces a bijection between the set \[ \{x\in \fa\cap \cA\mid \, \ord_v(x) =0 \,\,\,\forall \, v\in S, v\nmid \infty\} \] and the set $\sI_{\fA^{-1}}^S$ of integral ideals $\fb$ relatively prime to $S$ and  contained in $\fA^{-1}$. Define the Dirichlet series 
\begin{equation*}
\widetilde{\zeta}_S(\fa, s) \, = \, \Norm(\fa)^{s}[\sigma_{\fa}]\sum_{x\in \fa\cap \cA\cap U_S} [\rec_S(x)]\Norm(x)^{-s} \, = \,  \sum_{\fb\in \sI_{\fA^{-1}}^S} [\sigma_{\fb}^{-1}] \Norm(\fb)^{-s}.
\end{equation*}
For $L(1_{\cA}, \Phi_{\fa}, s)$ we obtain
\begin{equation*}
\Norm(\fa)^{s}L(1_{\cA}, \Phi_{\fa}, s) \, =\,\sum_{T'\subseteq T} (-1)^{\sharp(T')} \left( \prod_{\fq\in T'}\Norm(\fq)^{1-s} [\sigma_{\fq}^{-1}]\right)\widetilde{\zeta}_S\!\!\left(\fa\prod_{\fq\in T'}\fq\, s\right).
\end{equation*}
For a fixed subset $T'$ of $T$ the fractional ideals $\fa_1\prod_{\fq\in T'}\fq, \ldots, \fa_h\prod_{\fq\in T'}\fq$ are representatives of the different classes in $\Cl^+(F)$ as well, so we get
\begin{equation*}
\sum_{i=1}^h\widetilde{\zeta}_S\!\!\left(\fa_i\prod_{\fq\in T'}\fq, s\right) \, =\, \sum_{\fb}  [\sigma_{\fb}^{-1}] \Norm(\fb)^{-s}\, =\, \sum_{\sigma\in G}\zeta_S(\sigma,s)[\sigma^{-1}]
\end{equation*}
where in the second sum $\fb$ is taken over all integral ideals of $\cO_F$ that are relatively prime to $S$. We deduce
\begin{eqnarray*}
\zeta_{S, T}(K/F, s)&= & \,\sum_{T'\subseteq T} (-1)^{\sharp(T')} \left( \prod_{\fq\in T'}\Norm(\fq)^{1-s} [\sigma_{\fq}^{-1}]\right)\sum_{\sigma\in G}\zeta_S(\sigma,s)[\sigma^{-1}]\\
 & = & \sum_{i=1}^h\Norm(\fa_i)^s L(1_{\cA}, \Phi_{\fa_i}, s) 
\end{eqnarray*}
as desired.
\enddemo

{\em Proof of Prop.\ \ref{prop:zetaelement}.} As explained in Remark \ref{remark:stickelberger} (a) we can choose $S_3 =\emptyset$, $S_4=S_{\infty}$ 
in the definition of $\Theta_{S,T}(K/F,k)$. 

Choose ideles $a_1, \ldots, a_h\in (\bA_F^{\barT,\infty})^*$ whose components at places in $S$ are $=1$ and whose associated fractional ideals are $\fa_1, \ldots, \fa_h$. Then $\cF=\bigcup_{i=1}^h a_i U^{\barT,\infty}\subseteq (\bA_F^{\barT,\infty})^*/U^{\barT,\infty}$ is a fundamental domain for the action of $F^{\barT}/E_+$. Since $\cC(S, \bZ[G]) \cong \, \Coind^{F^*}_{F^{\barT}} \cC(S, \bZ[G])^{\barT,\infty}$ we have 
\begin{equation}
\label{delta10}
H^0(F^*, \cC(S, \bZ[G]))\, \cong \, H^0(F_+^{\barT}, \cC(S, \bZ[G])^{\barT,\infty}).
\end{equation}
Consider the pairing $\cC(S, \bZ[G])^{\barT,\infty}\times \cC_c(S, \bZ)^{\barT,\infty}\to  \cC_c(S, \bZ[G])^{\barT,\infty}$.   Passing to $E_+$-invariants and evaluating at $1_{\cF}\in H^0(E_+, \cC_c(S, \bZ)^{\barT,\infty})$ yields a map
\begin{equation}
\label{delta12}
H^0(E_+, \cC(S, \bZ[G])^{\barT,\infty}) \to H^0(E_+, \cC_c(S, \bZ[G])^{\barT,\infty}).
\end{equation}
The image of $\rec_{K/F}$ under the composition of \eqref{delta10}, \eqref{delta12} and \eqref{delta4a}
\begin{equation}
\label{delta12a}
H^0(F^*, \cC(S, \bZ[G]))\, \lra \, \cP_c(V, \bZ[G])^{\sL}
\end{equation}
is easily computed to be 
\begin{equation*}
\Phi_{K/F, k} \, =  \,\sum_{i=1}^h\, \Norm(\fa_i)^k \phi_{\fa_i} \otimes \Norm^{-k}.
\end{equation*}
Using the compatibility of the cap-product with restrictions and corestrictions one checks that the diagram 
\[
\begin{CD}
H^0(F^*, \cC(S, \bZ[G])) @> \partial^k >> H_{n-1}(F^{\barT}, \cP_c(V, \bZ[G])^{\sL})\\
@VV  \eqref{delta12a} V @AA \cor A\\
H^0(E_+, \cP_c(V, \bZ[G])^{\sL}) @> \cap \eta >> H_{n-1}(E_+, \cP_c(V, \bZ[G])^{\sL})
\end{CD}
\]
commutes. In particular we have $\partial^k(\rec_{K/F}) = \cor(\Phi_{K/F, k} \cap \eta)$. 

The specialization at $Q$ map $\spez\colon \Hom(\sR, \cN^{\sL}) \to \cN^{\sL}, f\mapsto f(Q)$ together with the inclusion $E_+\hookrightarrow \Gamma^{\sL}$ induces a homomorphism
\begin{equation*}
\spez_*: H^{n-1}(\Gamma^{\sL}, \Hom(\sR,\cN^{\sL})) \, \lra \, H^{n-1}(E_+, \cN^{\sL}).
\end{equation*}
Moreover \eqref{solhupol2} induces a pairing 
\begin{equation*}
\cap: H^{n-1}(E_+, \Hom(\cP_c(V, \bZ)^{\sL}, \bZ)) \times H_{n-1}(E_+, \cN^{\sL})\, \lra\, R
 \end{equation*}
and it is easy to see that 
\begin{equation*}
\Eis_F \cap \cor(z)\, =\, \spez_*(\Xi_{\Sh,\sL})\cap z
\end{equation*}
for all $z\in H_{n-1}(E_+, \cP_c(V, R)^{\sL})$. It follows that
\begin{eqnarray*}
&& \Eis_F  \cap \partial^k(\rec_{K/F})\, =\,\Eis_F  \cap \cor(\Phi_{K/F, k} \cap \eta)\\
&& = \,(\spez_*(\Xi_{\Sh,\sL}) \cup \Phi_{K/F, k}) \cap \eta \, = \,  \Phi_{K/F, k} \cap (\spez_*(\Xi_{\Sh,\sL}) \cap \eta)
 \end{eqnarray*}
By (\cite{cdg}, Theorem 1.5) there exists a Shintani decomposition $\cA$ such that $1_{\cA}\in \cK^{\sL}$ represents the class $\eta \cap \ev_*(\Xi_{\Sh,\sL})\in H_0(E_+, \cN^{\sL})$. The assertion follows now from Lemmas \ref{lemma:shintani} and \ref{lemma:shintani2}.\enddemo

\subsection{Gross' tower of fields conjecture}
\label{subsection:grosstower}

As before let $K/F$ be a finite abelian extension with Galois group $G$ and let $S$ be a finite set of nonarchimedean places of $F$ containing all finite places that are ramified in $K$. For any place $v$ of $F$ we denote by $G_v\subseteq G$ the decomposition group at $v$. For $v\in S$ put $I_v =\Ker(\bZ[G]\to \bZ[G/G_v])$. Also for $v\in S_{\infty}$ let $\sigma_v$ be a generator of $G_v$ and define ideals $I_v^{\pm}\subseteq \bZ[G]$ as follows
\begin{equation*}
I_v^\pm  \,\,\, = \,\,\,  \langle [\sigma_v] \mp 1\rangle 
\end{equation*}
(so in particular $I_v^+ = 0$ and $I_v^- = 2\bZ[G]$ if the places of $K$ above $v$ are real).

\begin{theo}
\label{theo:highervanishing}
Let $v_0$ be a fixed archimedean place of $F$

\noi (a) Assume that (A1), (A2) and (A3a) hold. Then we have 
\[
\Theta_{S,T}(K/F, 0)\in \prod_{v\in S} I_v \cdot \prod_{v\mid\infty, v\ne v_0} I_v^+, \qquad 
2\Theta_{S,T}(K/F, 0)\in \prod_{v\in S} I_v \cdot \prod_{v\mid\infty} I_v^+.
\]
\noi (b) Let $k\in \bZ_{<0}$ and assume that (A1), (A2) and (A3) hold. If $k$ is even (resp.\ odd) then we have 
\begin{align*}
\Theta_{S,T}(K/F, k) & \in \prod_{v\mid\infty, v\ne v_0} I_v^+ \qquad (\mbox{resp. }  \Theta_{S,T}(K/F, k)\in \prod_{v\mid\infty, v\ne v_0} I_v^-
), \\
2\Theta_{S,T}(K/F, k) & \in \prod_{v\mid\infty} I_v^+ \qquad (\mbox{resp. }  \Theta_{S,T}(K/F, k)\in \prod_{v\mid\infty} I_v^-).
\end{align*}
\end{theo}

{\em Proof.} (a) Put $S_1 = S$, $S_2 =\emptyset$, $S_3 = S_{\infty}-\{v_0\}$ and $S_4 = \{v_0\}$. Let $\rho_{K/F}$ be the image $\rec_{K/F}$ under the composition 
\begin{eqnarray}
\label{higherordervanishing1a}
&&  H^0(F^*, \cC(S, \bZ[G])) \, \stackrel{\eqref{delta2}}{\lra} \, H_{n-1}(F^*, \cC_c(S, \bZ[G]))\\ 
&& \hspace{1cm} \,\stackrel{\eqref{augep3}_*}{\lra} \, H_{n-1}(F^*, \cC_c(S, \bZ[G])^{S_3}(\sgn_{S_3})) \nonumber\\
&& \hspace{1cm} \, \stackrel{\eqref{extensionsbyzero}_*}{\lra}\, 
H_{n-1}(F^*, \cC_c(S_1 , S_2, \bZ[G])^{S_3}(\sgn_{S_3})).\nonumber
\end{eqnarray}
Choose $Q\in \sF_{v_0}$ and let 
\begin{equation*}
\Delta_*\colon H_{n-1}(F^*, \cC_c(S_1 , S_2, \bZ[G])^{S_3}(\sgn_{S_3}))\to  H_{n-1}(F^{\barT\cup S_4}, \Ccd(V_{\bhatZ^T},  \bZ[G])(\sgn))
\end{equation*}
be the map \eqref{delta6b} (for $A=\bZ[G]$). Then we have 
\begin{equation*}
\Theta_{S, T}(K/F, 0)\,\, = \,\, \Eis_{F, T, Q}^0\cap \Delta_*(\rho_{K/F}).
\end{equation*}
Let $R= \bZ[G]/\prod_{v\in S} I_v \cdot \prod_{v\mid\infty, v\ne v_0} I_v^+$ and let 
$\pi\colon \bZ[G]\to R$ be the projection. Since the local components $(\wcdot, K/F)_v$ of $\rec_{K/F}$ are trivial modulo $I_v$ (resp.\ modulo $I_v^+$) for all $v\in S$ (resp.\ $v\in S_{\infty}$) we can apply Prop.\ \ref{prop:higherordervanishing2}, so we get that $\rho_{K/F}$ is mapped to $0$ under the induced map
\begin{equation}
\label{higherordervanishing1b}
\pi_*\colon H_{n-1}(F^*, \cC_c(S_1 , S_2, \bZ[G])^{S_3}(\sgn_{S_3})) \lra
 H_{n-1}(F^*, \cC_c(S_1 , S_2, R)^{S_3}(\sgn_{S_3})).
\end{equation}
Hence we obtain
\begin{equation*}
\pi(\Theta_{S, T}(K/F, 0))\, = \, \Eis_{F, T, Q}\cap \Delta_*(\pi_*(\rho_{K/F})) \, =\, 0
\end{equation*}
i.e.\ $\Theta_{S, T}(K/F, 0)\in \prod_{v\in S} I_v \cdot \prod_{v\mid\infty, v\ne v_0} I_v^+$. 
 
 For the second assertion put $S_3= S_{\infty}$, $S_4=\emptyset$ and let $v$ be any archimedan place. We now work with the class $\Eis_{F, T, v}^0$ (see (\ref{e:defeisv}); note that it is a cohomology class for the group $F^{\barT}$ and not just $F^{\barT,v}$). We have 
\begin{equation*}
2 \, \Theta_{S, T}(K/F, 0)\,\, = \,\, \Eis_{F, T, v}^0\cap \Delta_*(\rho_{K/F}).
\end{equation*}
and therefore we can conclude as above that $2 \, \Theta_{S, T}(K/F, 0)\in I := \prod_{w \in S} I_w \prod_{w \mid \infty} I_w^{+}$ since $\rho_{K/F}\equiv 0 \pmod I$.
 
The proof of (b) is similar. Let $S_3$, $S_4$ be as above but we now put $S_1= \emptyset$, $S_2=S$. Moreover we set $R= \bZ[G]/\prod_{v\mid\infty, v\ne v_0} I_v^+$ if $k$ is even and $R= \bZ[G]/\prod_{v\mid\infty, v\ne v_0} I_v^-$ if $k$ is odd and we let $\pi: \bZ[G]\to R$ be the projection. Again $\rho_{K/F}$ denotes the image $\rec_{K/F}$ under \eqref{higherordervanishing1a}. Since $(-1, K/F)_v = \sigma_v \equiv 1 = -\sgn_v(-1)^{k+1}$ modulo $I_v^+$ if $k$ is even (resp.\ $(-1, K/F)_v \equiv -\sgn_v(-1)^{k+1}$ modulo $I_v^-$ if $k$ is odd) by applying Prop.\ \ref{prop:higherordervanishing2} we see that $\rho_{K/F}$ is mapped to $0$ under \eqref{higherordervanishing1b}. Since the map \eqref{delta6} factors through \eqref{higherordervanishing1a} we conclude
\begin{equation*}
\pi(\Theta_{S, T}(K/F, k))\,\, = \,\, \Eis_{F, T, Q}\cap \partial^k(\pi_*(\rec_{K/F}))  \,=\, 0
\end{equation*}
hence $\Theta_{S, T}(K/F, k)\in \prod_{v\mid\infty, v\ne v_0} I_v^+$.\enddemo

Assume now that $L\supseteq K\supseteq F$ is a tower of finite abelian extensions (i.e\ $L/F$ is abelian). Let $G$, $\widetilde{G}$ and $H$ be the Galois groups of $K/F$, $L/F$ and $L/K$ respectively and let $I$ be the kernel of the canonical projection $\bZ[\widetilde{G}] \to \bZ[G]$. We assume that $L/F$ is unramified outside $S$ and denote by $r$ be the number of places in $S$ that split completely in $K$. Furthermore  let $s$ be the number of archimedean places in $F$ that   split completely in $K$. Since $I_v\subseteq I$ (resp.\ $I_v^+ \subseteq I$) if $v\in S$ (resp.\ $v\in S_{\infty}$) splits completely in $K$ we obtain in particular

\begin{corollary}
\label{corollary:grosstower}
(a) Assume that (A1), (A2) and (A3a) hold. Then 
\[
\Theta_{S,T}(L/F, 0)\in I^{r+ \min(s, n-1)} \qquad \text{ and } \qquad 2\Theta_{S,T}(L/F, 0)\in I^{r+s}.
\]
\noi (b) Let $k\in \bZ_{<0}$ be even and assume that (A1), (A2) and (A3) hold. Then 
\[
\Theta_{S,T}(L/F, k)\in I^{\min(s, n-1)} \qquad \text{ and } \qquad 2\Theta_{S,T}(L/F, k)\in I^{s}.
\]
\end{corollary}

\begin{remark}
\label{remarks:highervanishing1}
\rm In \cite{gross2} Gross  conjectured that $\Theta_{S,T}(L/F, 0)\in I^m$ with \[ m= \min(r+s, \sharp(S\cup S_{\infty})-1). \] In fact if $r\ge 2$ then (\cite{gross2}, Conj.\ 7.6) predicts $\Theta_{S,T}(L/F, 0)\in I^2$ and the strengthening $\Theta_{S,T}(L/F, 0)\in I^{m}$ was formulated as a ``guess'' in \cite{gross2} on top of page 195. Therefore our result implies Gross' conjecture unless $K$ is totally real and not all the places of $S$ split completely in $K$; in this case, our exponent is one less than that which Gross predicts for the ``2-part" of the conjecture.
\end{remark}

\section{Conjectural construction of Gross-Stark units}
\label{section:grossstark}

In this section we give a cohomological interpretation of the {\it Gross-Stark units} introduced in \cite{dasgupta1}. 

Let $F$ again denote a totally real number field of degree $n$ over $\bQ$, let $K/F$ be a finite abelian extension with Galois group $G$ and let $S$ denote a set of nonarchimedean places of $F$ containing all places that are ramified in $K$. Let $\fp \not \in S$ be a nonarchimedean place of $F$ that splits completely in $K$. We let $T$, $\barT$ be as in the the beginning of section \ref{subsection:totalreal} such that $\fp\not\in \barT$ and such that $T$ satisfies assumptions (A1), (A2) and (A3a). Since the local norm residue symbol for $K/F$ at $\fp$ is trivial we omit it from the reciprocity map, i.e.\ we consider the homomorphism
\begin{equation*}
\rec_{K/F}^{\fp}\colon (\bA_F^{\fp})^*\,\lra\, \bZ[G]^*, \qquad x=(x_v)_{v\ne \fp} \mapsto \prod_{v\ne \fp} (x, K/F)_v.
\end{equation*}
We view it as an element of $H^0(F^*, \cC(S, \bZ[G])^{\fp})$ and denote by 
\begin{equation*}
\rho_{K/F} \in H_n(F^*,\cC_c(S, \bZ[G])^{\fp})
\end{equation*}
its image under \eqref{deltaS2}. 

For a locally profinite abelian group $A$ the bilinear map $\otimes\colon A \times \bZ[G] \to A\otimes \bZ[G]$ induces a bilinear map 
\begin{equation*}
\Ccd(F_{\fp}, A)\otimes\, \cC_c(S, \bZ[G])^{\fp}\,\lra \, \cCcd(\{\fp\}, S, A\otimes \bZ[G])
\end{equation*}
(compare \eqref{multfunctions4}) hence induces a pairing 
\begin{equation*}
H^i(F^*, \Ccd(F_{\fp}, A)) \times H_j(F^*,\cC_c(S, \bZ[G])^{\fp})\to H_{j-i}(F^*,\cCcd(\{\fp\}, S, A\otimes \bZ[G])).
\end{equation*}
In particular we can consider
\begin{equation*}
c_{\univ}\cap \rho_{K/F}\in H_{n-1}(F^*, \cCcd(\{\fp\},S, \bZ[G]\otimes F_{\fp}^*))
\end{equation*}
where $c_{\univ} = c_{\id} \in H^1(F^*, \Ccd(F_{\fp}, F_{\fp}^*))$ is the class \eqref{cyclehom} attached to the identity 
$\id\colon F_{\fp}^*\to F_{\fp}^*$.
Choose $v\in S_{\infty}$ and $Q\in \sF_v$. Applying the map induced by \eqref{delta6b} we get a class 
\begin{equation*}
\Delta_*(c_{\univ}\cap \rho_{K/F})\in H_{n-1}(F^{\barT,v}, \Ccd(V_{\bhatZ^T}, F_{\fp}^*\otimes\bZ[G])(\sgn)).
\end{equation*}

Now the canonical pairing (set $A = F_\fp^* \otimes \bZ[G]$)
\begin{equation*}
\Hom(C_c(V_{\bhatZ^T}, \bZ), \bZ)\times \Ccd(V_{\bhatZ^T}, A)\lra A, \qquad (\mu, f) \mapsto \mu_A(f)
\end{equation*}
induces via cap-product a pairing
\begin{equation}
\label{capqmeas}
\cap: H^{n-1}(F^{\barT,v}, \Hom(C_c(V_{\bhatZ^T}, \bZ), \bZ)(\sgn))\times H_{n-1}(F^{\barT,v}, \Ccd(V_{\bhatZ^T}, A)(\sgn))\to A.
\end{equation}
Applying \eqref{capqmeas} with the Eisenstein cocycle $\Eis_F^0=\Eis_{F, T, Q}^0$
and $\Delta_*(c_{\univ}\cap \rho_{K/F})$
 we obtain an element $u_{K/F} =  u_{K/F, S, T}\in F_{\fp}^*\otimes \bZ[G]$:
\begin{equation}
\label{grossstarklocal}
u_{K/F, S,T} \,= \, \sum_{\sigma\in G} u_{S,T, \sigma} \otimes [\sigma^{-1}]  \, =\, \Eis_F^0 \cap \Delta_*(c_{\univ}\cap \rho_{K/F}).
\end{equation}

Let $\sE_{\fp}$ denote the subgroup 
\begin{equation*}
\sE_{\fp} \, =\, \{x\in K^*\mid |x|_w=1 \, \mbox{for all places $w$ of $H$ with $w\nmid \fp$}\} 
\end{equation*}
of $K^*$ and choose a prime $\fP$ of $K$ above $\fp$. Since $F_{\fp} = K_{\fP}$ we can view $K$ as a subfield of $F_{\fp}$. The following conjecture is a strengthening of the Brumer--Stark conjecture, as formulated by Tate, and a conjecture of Gross (see \cite{tate}, \cite{gross2}; compare also \cite{dasgupta1}, section 2).

\begin{conjecture}
\label{conjecture:grossstark}
\noi (i) For all $\sigma\in G$ we have $u_{S, T,\sigma}\in \sE_{\fp}$.

\noi (ii) We have $u_{S, T,\sigma}\equiv 1\mod T$, i.e.\ for all $\sigma\in G$ and all $\fq\in T$ we have $u_{S, T,\sigma}\equiv 1\mod \fq\cO_H$.

\noi (iii) For all $\sigma, \tau\in G$ we have $\tau(u_{S, T,\sigma}) = u_{S, T,\tau\sigma}$.
\end{conjecture}

\begin{remark} \rm Conditions (i) and (iii) may be rephrased by conjecturing that $u_{K/F, S, T}$ lies in the image of the diagonal embedding 
\begin{equation*}
\sE_{\fp} \,\lra \, \prod_{\sigma\in G} K_{\sigma\fP}^* = \prod_{\sigma\in G} F_{\fp}^*.
\end{equation*}
Hence Conjecture \ref{conjecture:grossstark} is independent of the choice of $\fP$.
\end{remark}

\begin{prop}
\label{prop:brumerstark}
(a) For $\sigma\in G$ we have $\ord_{\fp}( u_{S, T,\sigma}) \,\, = \,\, \zeta_{S,T}(\sigma,0)$.

\noi (b) Let $L/F$ be an abelian extension with $L\supseteq K$ and put $\widetilde{G} =\Gal(L/F)$. Assume that $L/F$ is unramified outside $S$ and that $\fp$ splits completely in $L$. Then we have  
\[
u_{S, T,\sigma}= \prod_{\tau\in \widetilde{G}, \tau|_K = \sigma}u_{S, T,\tau}.
\]
\noi (c) Let $\fq$ be a nonarchimedean place of $F$ with $\fq\not\in S\cup \barT\cup \{\fp\}$ and put $S'=S\cup \{\fq\}$. Then we have $u_{S', T,\sigma} \,\, = u_{S, T,\sigma}u_{S, T,\sigma_{\fq}\sigma}^{-1}$.

\noi (d) Assume that $K$ has a real archimedean place $w$ with $w\nmid v$. Then $u_{S, T,\sigma}=1$ for all $\sigma\in G$.  

\noi (e) Let $L/F$ be an abelian extension with Galois group $\widetilde{G}$ and assume $K\subseteq L$ and that $L/F$ is unramified outside $S'= S\cup\{\fp\}$. Let \[ \rec_{\fp}= (\wcdot, L/F)_{\fp}\colon F_{\fp}^*\to \Gal(L/K)\subseteq \widetilde{G} \] be the $\fp$-component of the reciprocity map. Then we have \[ \rec_{\fp}( u_{S, T,\sigma}) \,\, = \,\, \prod_{\tau\in \widetilde{G}, \tau|_H = \sigma} \tau^{\zeta_{S',T}(\tau,0)}. \]
\end{prop}

{\em Proof.} (a) This follows from Lemma \ref{lemma:capord} and Prop.\ \ref{prop:zetaelement}. Indeed, using Lemma \ref{lemma:capord} it is easy to see that the image of $u_{K/F, S,T}$ under the map $\ord_{\fp}\otimes \id_{\bZ[G]} \colon F_{\fp}^*\otimes \bZ[G]\to \bZ[G]$ is $=\Theta_{S, T}(K/F, 0)$, i.e.\ the element we obtain by replacing $c_{\univ}$ by $c_{\ord_{\fp}}$ in \eqref{grossstarklocal} above.

\noi (b) is obvious and (c) can be deduced from the commutativity of diagram \eqref{changeS2} in the same way as the formula \eqref{changeS0} in Remark \ref{remark:stickelberger}(c).

\noi (d) Put $v' = v|_F$. It is easy to see that the map 
\begin{eqnarray*}
&& H_n(F^*,\cC_c(S, \bZ[G])^{\fp})\to H_{n-1}(F^{\barT,v}, \Ccd(V_{\bhatZ^T}, F_{\fp}^*\otimes\bZ[G])(\sgn)),\\
&& \hspace{3cm} x\,\mapsto \, \Delta_*(c_{\univ}\cap x)
\end{eqnarray*}
factors through the right vertical map of the diagram 
\begin{equation*}
\begin{CD}
H^0(F^*, \cCd(S, \bZ[G])^{\fp}) @>\eqref{deltaS2}>> H_n(F^*, \cCcd(S,  \bZ[G])^{\fp})\\
@VV \eqref{augep2}_* V @VV\eqref{augep3}_* V\\
H^0(F^*, \cCd(S, \bZ[G])^{\fp, v'}(\sgn_{v'})) @ > \eqref{delta2infty} >> H_n(F^*, \cCcd(S, \bZ[G])^{\fp, v'}(\sgn_{v'}))).
\end{CD}
\end{equation*}
Since $(\wcdot, K/F)_{v'} = 1$ the homomorphism $\rec_{K/F}^{\fp}$ lies in the kernel of the left vertical map. Consequently, we obtain $\Delta_*(c_{\univ}\cap \rho_{K/F})=0$, hence $u_{K/F, S, T}$ is trivial.

\noi (e) Let $I$ be the kernel of the canonical projection $\pi\colon\bZ[\widetilde{G}]\to \bZ[G]$ and let $H=\Gal(L/K)$. The homomorphism
\begin{equation*}
\bZ[\widetilde{G}]\, \lra \, \widetilde{G}\otimes \bZ[G], \qquad \sum_{\tau\in \widetilde{G}}n_{\tau} [\tau]\mapsto \sum_{\sigma\in G} \left(\prod_{\tau\in \widetilde{G},\tau|_H=\sigma} \tau^{n_{\tau}}\right) \otimes [\sigma]
\end{equation*}
maps $I$ into $H\otimes \bZ[G]\subseteq \widetilde{G}\otimes \bZ[G]$ and induces an isomorphism
\begin{equation}
\label{gross-stark5}
I/I^2 \, \lra \, H\otimes \bZ[G].
\end{equation}
Thus the assertion will follow once we show that the image of $u_{K/F, S,T}$ under \[ \rec_{\fp}\otimes \id_{\bZ[G]}\colon F_{\fp}^*\otimes \bZ[G] \longrightarrow H\otimes \bZ[G] \] is equal to the image of $\Theta_{L/F, S',T}\mod I^2$ under the isomorphism \eqref{gross-stark5}. Firstly, note  
\begin{equation*}
(\rec_{\fp}\otimes \id_{\bZ[G]})(u_{K/F, S,T}) \,=\, \Eis_F^0 \cap \Delta_*(c_{\rec_{\fp}}\cap \rho_{K/F}).
\end{equation*}
In order to compute the right hand side further we apply Prop.\ \ref{prop:higherordervanishing} to the ring $R= \bZ[\widetilde{G}]/I^2$, the ideal $\fa = I/I^2$ and the character \[ \chi: \bA_F^*/U^{S'}F^*\lra R^*, \qquad x\mapsto \rec_{L/F}(x)\mod I^2. \] Note that $\chi$ mod $\fa = \rec_{K/F}$ so with the notation as in Prop.\ \ref{prop:higherordervanishing} we have $\kappa_{\chi, \{\fp\}} =  c_{\der\!\chi_{\fp}}\cap \rho_{K/F}$. Consider the commutative diagram 
\[
\begin{CD}
H\times \bZ[G] @> \otimes >> H\otimes \bZ[G]\\
@VVV @AA \eqref{gross-stark5} A\\
I/I^2 \times R/I @>\mult >> I/I^2
\end{CD}
\]
where the first vertical arrow is given on the first factor by $H\to I/I^2, \tau\mapsto [\tau] - 1\mod I^2$ and by $R/I\cong \bZ[G]$ on the second factor. Note that the image of 
$c_{\rec_{\fp}}$ under the map $H^1(F_{\fp}^*, C_c(F_{\fp}, H))\to H^1(F_{\fp}^*, C_c(F_{\fp}, I/I^2))$ induced by $H\to I/I^2$ as above is $c_{\der\!\chi_{\fp}}$. The functorial properties of the cap- and cup-product yield $\kappa_{\chi, \{\fp\}} = c_{\rec_{\fp}}\cap \rho_{K/F}$. Applying $\Delta^*$ to both sides and taking the cap-product with $\Eis_F^0$ completes the proof.\enddemo 

\begin{remarks}
\label{remarks:gross-stark}
\rm (a) As mentioned in the introduction, elements $u_{S,T,\sigma}\in F_{\fp}^*$ for which the Conjecture \ref{conjecture:grossstark} is expected to hold have been introduced previously in (\cite{dasgupta1}, Def.\ 3.18). The definition given there is a rather elaborate construction in terms of a certain $p$-adic measure involving Shintani zeta functions (it does not make use of group homology and cohomology). One can show by a tedious but ultimately straightforward computation that our cohomological construction coincides with the one in loc.\ cit.\ 

\noi (b) Of course we expect $u_{S,T, \sigma}$ to be independent of the choice of $v$ and $Q$. In particular we expect that $u_{S,T, \sigma}=1$ for all $\sigma \in G$  if $K$ has a real  place; indeed, this follows from  Conjecture \ref{conjecture:grossstark}(i) since $\sE_{\fp} = 1$ unless $K/F$ is totally imaginary.

\noi (c) Following the arguments of \cite{dasgupta1}, pp.\ 262--264 one can show that the properties of Prop.\ \ref{prop:brumerstark} determine the elements $u_{S,T, \sigma}$ ``almost'' uniquely. More precisely assume that $K/F$ is totally imaginary and choose an archimedean place $v'$ of $F$ with $v'\ne v$. Let $\tau\in G_{v'}\subseteq G$ denote the complex conjugation associated to $v'$. If $\fq_1, \ldots, \fq_m$ are nonarchimedan places of $F$ with $\fq_i\ne S\cup \barT$ and $\sigma_{\fq_i} =\tau$ for all $i=1,\ldots, m$ then by Prop.\ \ref{prop:brumerstark} (b), (c), (d) we have for $S'= S\cup \{\fq_1, \ldots, \fq_m\}$
\[
u_{S',T, \sigma} \,= \, u_{S,T, \sigma}^{2^m}.
\]
As in \cite{dasgupta1} it follows that by enlarging $S$ further and further with primes $\fq$ as above (using Lemma 5.17 of loc.\ cit.) the element $u_{S,T, \sigma}$ is uniquely determined up to a root of unity by the properties (a) and (e) of Prop.\ \ref{prop:brumerstark}.

\noi (d) Assume that $\fp$ lies over a prime number $p$ such that $S$ contains all other primes of $F$ above $p$. It is easy to show (see e.g.\ the forthcoming work \cite{dasspi}) that we have
\begin{equation}
\label{formalgross}
\log_p(\Norm_{F_{\fp}/\bQ_p}(u_{S,T,\sigma})) \, =\, \zeta_{S, T,p}'(\sigma, 0)
\end{equation}
where $\zeta_{S,T,p}(\sigma, s)$, $s\in \bZ_p$ is the partial $p$-adic zeta function (it satisfies the interpolation property $\zeta_{S,T,p}(\sigma, n) = \zeta_{S,T}(\sigma, n)$ for all $n\in \bZ_{\le 0}$ such that $n \equiv 0 \pmod{p-1}$). Assume now that $p\ge 3$, $F_{\fp} = \bQ_p$ and that $\fp$ is the only prime in $S_p$ that splits completely in $K$. Combining \eqref{formalgross} with the main results of \cite{ddp} and \cite{ventullo} we obtain that at least some power of $u_{S,T,\sigma}$ is contained in $K$.
\end{remarks}

 \bigskip

 S.~D. : \textsc{Dept. of Mathematics, University of California  Santa Cruz, USA} 

 \textit{E-mail} \textbf{sdasgup2 (at) ucsc (dot) edu}
 
 \bigskip
 
 M.~S. : \textsc{Fakult{\"a}t f{\"u}r Mathematik, Universit{\"a}t Bielefeld, Germany} 

 \textit{E-mail} \textbf{mspiess (at) math (dot) uni-bielefeld (dot) de}
 
\end{document}